\documentclass[reqno]{amsart}%
\usepackage{amsmath}
\usepackage{graphicx}
\usepackage{array}
\usepackage{epsfig}
\usepackage{float,array}
\usepackage{amsfonts}
\usepackage{amssymb}
\usepackage[dvipdfm]{hyperref}
\usepackage[all]{xy}%
\providecommand{\U}[1]{\protect\rule{.1in}{.1in}}
\newtheorem{theorem}{Theorem}[section]

\newtheorem{lemma}[theorem]{Lemma}
\newtheorem{proposition}[theorem]{Proposition}

\theoremstyle{definition}
\newtheorem{definition}[theorem]{Definition}

\newtheorem{remark}[theorem]{Remark}
\newtheorem{note}[theorem]{Note}
\numberwithin{equation}{section}
\begin{document}
\author[D.I. Dais]{Dimitrios I. Dais}
\address{University of Crete, Department of Mathematics, Division Algebra and Geometry,
Knossos Avenue, P.O. Box 2208, GR-71409, Heraklion, Crete, Greece}
\email{ddais@math.uoc.gr}
\title[Toric log Del Pezzos with $\rho=1$ and index $\leq3$]{Classification of Toric log Del Pezzo Surfaces \\having Picard Number $1$ and Index $\leq3$}
\dedicatory{To Professor Friedrich Hirzebruch on the occasion of his eightieth birthday}\date{}
\maketitle

\begin{abstract}
Toric log Del Pezzo surfaces with Picard number $1$ have been completely
classified whenever their index is $\leq2.$ In this paper we extend the
classification for those having index $3.$ We prove that, up to isomorphism,
there are exactly $18$ surfaces of this kind.

\end{abstract}

\section{Introduction\label{INTRO}}

\noindent{}Smooth compact toric surfaces belong to the basics in the framework
of toric geometry. They are rational surfaces (i.e., of Kodaira dimension
$-\infty$) defined by $2$-dimensional complete fans which are composed of
basic cones, and can therefore be studied by means of handy combinatorics (see
\cite[Theorem 1.28, pp. 42-43]{Oda}). Of course, unlike the smooth compact
complex surfaces having Kodaira dimension $\geq0,$ they do not possess
uniquely determined \textit{minimal models}. Nevertheless, the set of their
minimal models consists of the projective plane $\mathbb{P}_{\mathbb{C}}^{2}$
together with the \textit{Hirzebruch surfaces}%
\[
\mathbb{F}_{\kappa}:=\left\{  \left.  (\left[  z_{0}:z_{1}:z_{2}\right]
,\left[  t_{1}:t_{2}\right]  )\in\mathbb{P}_{\mathbb{C}}^{2}\times
\mathbb{P}_{\mathbb{C}}^{1}\ \right\vert \ z_{1}t_{1}^{\kappa}=z_{2}%
t_{2}^{\kappa}\right\}  ,\ \ \kappa\in\mathbb{Z}_{\geq0},
\]
for $\kappa$ $\neq1$ (cf. \cite{Hirzebruch1}, \cite[\S 2.5]{Fulton},
\cite[\S 1.7]{Oda}), and it is known how one can pass from one minimal model
to another by a finite succession of \textit{elementary transformations}.

In contrast to this classical point of view, taking into account the fact that
the \textit{anti-Kodaira dimension} of smooth compact toric surfaces is $2,$
and switching to the so-called \textit{antiminimal} and \textit{anticanonical
models} (in the sense of Sakai \cite[\S 7]{Sakai2}, \cite[Appendix]{Sakai3}),
one obtains surfaces which are \textit{uniquely determined} up to isomorphism.
However, since these models are mostly \textit{singular}, in order to follow
this choice we need a more systematic study of \textit{singular} compact toric surfaces.

A graph-theoretic method of classifying (not necessarily smooth) compact toric
surfaces \textit{up to isomorphism} (generalizing Oda's graphs \cite[pp.
44-46]{Oda}) has been proposed in \cite[\S 5]{Dais}: Two compact toric
surfaces are isomorphic to each other if and only if their vertex singly- and
edge doubly-weighted circular graphs (\textsc{wve}$^{2}$\textsc{c}%
-\textit{graphs}, for short) are isomorphic (see below Theorem
\ref{CLASSIFTHM}).

In addition, by \cite[Theorem 4.3, pp. 398-399]{Sakai1} the anticanonical
models of {}smooth compact toric surfaces have to be log Del Pezzo surfaces.
(A compact complex surface $X$ with at worst log terminal singularities, i.e.,
quotient singularities, is called \textit{log Del Pezzo} \textit{surface} if
its anticanonical divisor $-K_{X}$ is a $\mathbb{Q}$-Cartier ample divisor.
The \textit{index} of such a surface is defined to be the smallest positive
integer $\ell$ for which $\ell K_{X}$ is a Cartier divisor. The family of log
Del Pezzo surfaces of fixed\textit{ }index $\ell$ is known to be bounded (see
\cite[Theorem 2.1, p. 332]{Borisov})).

Consequently, it seems to be rather interesting to classify toric log Del
Pezzo surfaces of given index $\ell$ up to isomorphism. A first attempt to
understand the combinatorial complexity of \ this classification problem
includes naturally the investigation of the case in which the Picard number
$\rho\left(  X_{\Delta}\right)  :=$ rank$\left(  \text{Pic}\left(  X_{\Delta
}\right)  \right)  $ of surfaces $X_{\Delta}$ of this kind (associated to
complete fans $\Delta$ in $\mathbb{R}^{2}$) equals $1.$ In this case,
$X_{\Delta}$'s turn out to be weighted projective planes or quotients thereof
by a finite abelian group. Let us first recall what is known for indices
$\ell\leq2$:

\begin{theorem}
\label{THM1}Up to isomorphism, there are exactly $5$ toric log del Pezzo
surfaces with Picard number $1$ and index $\ell=1,$ namely
\setlength\extrarowheight{2pt}%
\[%
\begin{tabular}
[c]{|c||c|c|c|c|c|}\hline
\emph{No.} & \emph{(i)} & \emph{(ii)} & \emph{(iii)} & \emph{(iv)} &
\emph{(v)}\\\hline
$X_{\Delta}$ & $%
\begin{array}
[c]{c}%
\mathbb{P}_{\mathbb{C}}^{2}%
\end{array}
$ & $\mathbb{P}_{\mathbb{C}}^{2}/(\mathbb{Z}/3\mathbb{Z})$ & $\mathbb{P}%
_{\mathbb{C}}^{2}(1,1,2)$ & $\mathbb{P}_{\mathbb{C}}^{2}(1,1,2)/(\mathbb{Z}%
/2\mathbb{Z})$ & $\mathbb{P}_{\mathbb{C}}^{2}(1,2,3)$\\\hline
\end{tabular}
\ \ \ \ \ \ \ \ \ \ \ \
\]
\setlength\extrarowheight{-2pt}
whose \textsc{wve}$^{2}$\textsc{c}-graphs are illustrated in \cite[Figure 8,
p. 108]{Dais}.
\end{theorem}

\begin{theorem}
\label{THM2}Up to isomorphism, there are exactly $7$ toric log del Pezzo
surfaces with Picard number $1$ and index $\ell=2,$ namely
\setlength\extrarowheight{2pt}%
\[%
\begin{tabular}
[c]{|c|c|c|c|}\hline
\emph{No.} & $X_{\Delta}$ & \emph{No.} & $X_{\Delta}$\\\hline\hline
\emph{(i)} & $%
\begin{array}
[c]{c}%
\mathbb{P}_{\mathbb{C}}^{2}(1,1,4)
\end{array}
$ & \emph{(iv)} & $\mathbb{P}_{\mathbb{C}}^{2}(1,2,3)/(\mathbb{Z}%
/2\mathbb{Z})$\\\hline
\emph{(ii)} & $%
\begin{array}
[c]{c}%
\mathbb{P}_{\mathbb{C}}^{2}(1,4,5)
\end{array}
$ & \emph{(v)} & $\mathbb{P}_{\mathbb{C}}^{2}(1,1,2)/(\mathbb{Z}/4\mathbb{Z}%
)$\\\hline
\emph{(iii)} & $%
\begin{array}
[c]{c}%
\mathbb{P}_{\mathbb{C}}^{2}(1,3,8)
\end{array}
$ & \emph{(vi)} & $\mathbb{P}_{\mathbb{C}}^{2}(1,2,1)/(\mathbb{Z}%
/4\mathbb{Z})$\\\hline
&  & \emph{(vii)} & $%
\begin{array}
[c]{c}%
\mathbb{P}_{\mathbb{C}}^{2}(1,1,4)/(\mathbb{Z}/3\mathbb{Z})
\end{array}
$\\\hline
\end{tabular}
\ \ \ \ \ \ \ \ \ \ \
\]
\setlength\extrarowheight{-2pt}
whose \textsc{wve}$^{2}$\textsc{c}-graphs are illustrated in \cite[Figure 11,
p. 111]{Dais}.
\end{theorem}

\noindent{}\noindent{}In the present paper we extend these results also for
index $3$ by the following:

\begin{theorem}
\label{THM3}Up to isomorphism, there are exactly $18$ toric log del Pezzo
surfaces with Picard number $1$ and index $\ell=3,$ namely
\setlength\extrarowheight{2pt}%
\[%
\begin{tabular}
[c]{|c|c|c|c|}\hline
\emph{No.} & $X_{\Delta}$ & \emph{No.} & $X_{\Delta}$\\\hline\hline
\emph{(i)} & $\mathbb{P}_{\mathbb{C}}^{2}(1,1,3)$ & \emph{(x)} &
$\mathbb{P}_{\mathbb{C}}^{2}(1,5,9)$\\\hline
\emph{(ii)} & $\mathbb{P}_{\mathbb{C}}^{2}(1,3,4)$ & \emph{(xi)} &
$\mathbb{P}_{\mathbb{C}}^{2}(1,2,9)$\\\hline
\emph{(iii)} & $\mathbb{P}_{\mathbb{C}}^{2}(2,3,5)$ & \emph{(xii)} &
$\mathbb{P}_{\mathbb{C}}^{2}(1,2,3)/(\mathbb{Z}/3\mathbb{Z})$\\\hline
\emph{(iv)} & $\mathbb{P}_{\mathbb{C}}^{2}(1,1,2)/(\mathbb{Z}/3\mathbb{Z})$ &
\emph{(xiii)} & $\mathbb{P}_{\mathbb{C}}^{2}(1,1,2)/(\mathbb{Z}/2\mathbb{Z}%
)\times(\mathbb{Z}/3\mathbb{Z})$\\\hline
\emph{(v)} & $\mathbb{P}_{\mathbb{C}}^{2}(1,1,6)$ & \emph{(xiv)} &
$\mathbb{P}_{\mathbb{C}}^{2}(1,1,6)/(\mathbb{Z}/2\mathbb{Z})$\\\hline
\emph{(vi)} & $\mathbb{P}_{\mathbb{C}}^{2}(1,6,7)$ & \emph{(xv)} &
$\mathbb{P}_{\mathbb{C}}^{2}(1,4,15)$\\\hline
\emph{(vii)} & $\mathbb{P}_{\mathbb{C}}^{2}(1,3,4)/(\mathbb{Z}/2\mathbb{Z})$ &
\emph{(xvi)} & $\mathbb{P}_{\mathbb{C}}^{2}(1,1,3)/(\mathbb{Z}/5\mathbb{Z}%
)$\\\hline
\emph{(viii)} & $\mathbb{P}_{\mathbb{C}}^{2}(1,2,3)/(\mathbb{Z}/3\mathbb{Z})$
& \emph{(xvii)} & $\mathbb{P}_{\mathbb{C}}^{2}(1,2,9)/(\mathbb{Z}%
/2\mathbb{Z})$\\\hline
\emph{(ix)} & $\mathbb{P}_{\mathbb{C}}^{2}/(\mathbb{Z}/9\mathbb{Z}%
)$ & \emph{(xviii)} & $\mathbb{P}_{\mathbb{C}%
}^{2}(1,1,6)/(\mathbb{Z}/4\mathbb{Z})$\\\hline
\end{tabular}
\ \ \ \ \ \ \ \ \ \
\]
\setlength\extrarowheight{-2pt}
whose \textsc{wve}$^{2}$\textsc{c}-graphs are illustrated below in
\emph{Figure \ref{Fig.3}}.
\end{theorem}

\noindent{}The paper is organized as follows: In \S \ref{2DIMTORSING} we focus
on the properties of the two non-negative, relatively prime integers
$p=p_{\sigma}$ and $q=q_{\sigma}$ which parametrize the $2$-dimensional,
rational, strongly convex polyhedral cones $\sigma,$ and recall how they are
involved in Hirzebruch's minimal desingularization \cite{Hirzebruch2} of the
$2$-dimensional cyclic quotient singularities orb$\left(  \sigma\right)  \in$
Spec$(\mathbb{C}[\sigma^{\vee}\cap\mathbb{Z}^{2}])$ for $q>1.$ In section
\ref{LOCALINDICES} we give necessary and sufficient arithmetical conditions
for the local indices $l=l_{\sigma}$ to be $1$ or $3.$ Sections
\ref{COMPACTTS} and \ref{TLDPSURFACES} are devoted to a detailed description
of compact toric surfaces and of those which are log Del Pezzo surfaces. Some
key-lemmas of combinatorial nature concerning compact toric surfaces with
Picard number $1$ are presented in \S \ref{CTSWITHPIC1}. Based on the results
of \S \ref{LOCALINDICES}-\S \ref{CTSWITHPIC1} we explain how the
classification method works in \S \ref{STRATEGY}. The proof of Theorem
\ref{THM3} (which is somewhat longer than that of \ref{THM1} and \ref{THM2})
follows in four steps (in \S \ref{STEP1}-\S \ref{STEP4}). The first three
include the case by case determination of all \textquotedblleft
amissible\textquotedblright\ of triples of pairs $(p_{i},q_{i}),1\leq i\leq3,$
so that the induced toric log Del Pezzo surfaces $X_{\Delta}$ with Picard
number $\rho\left(  X_{\Delta}\right)  =1$ have index $\ell=3$. A minimal set
of pairwise non-isomorphic surfaces of this kind is sorted out in the fourth step.

We use tools only from the classical toric geometry, adopting the standard
terminology from \cite{Ewald}, \cite{Fulton}, and \cite{Oda} (and mostly the
notation introduced in \cite{Dais}).

\section{Two-dimensional toric singularities\label{2DIMTORSING}}

\noindent{}Let $\sigma=\mathbb{R}_{\geq0}\mathbf{n}+\mathbb{R}_{\geq
0}\mathbf{n}^{\prime}\subset\mathbb{R}^{2}$ be a 2-dimensional, rational,
strongly convex polyhedral cone. Without loss of generality we may assume that
$\mathbf{n}=\tbinom{a}{b}$, $\mathbf{n}^{\prime}=\tbinom{c}{d}\in
\mathbb{Z}^{2},$ and that both $\mathbf{n}$ and $\mathbf{n}^{\prime}$ are
primitive elements of $\mathbb{Z}^{2}$, i.e., gcd$\left(  a,b\right)  =1$ and
gcd$\left(  c,d\right)  =1.$

\begin{lemma}
\label{pequ}Consider $\kappa,\lambda\in\mathbb{Z},$ such that $\kappa
a-\lambda b=1.$ If $q:=\left\vert ad-bc\right\vert ,$ and $p$ is the unique
integer with
\[
0\leq p<q\text{ \ \ \ \emph{and} \ \ \ }\kappa c-\lambda d\text{ }\equiv
p\left(  \text{\emph{mod} }q\right)  ,
\]
then\emph{\ gcd}$\left(  p,q\right)  =1$, and there exists a primitive element
$\mathbf{n}^{\prime\prime}=\binom{e}{g}\in\mathbb{Z}^{2}$, such
that\ $\mathbf{n}^{\prime}=p\mathbf{n}+q\mathbf{n}^{\prime\prime}$ and
$\left\{  \mathbf{n},\mathbf{n}^{\prime\prime}\right\}  $ is a $\mathbb{Z}%
$-basis of $\mathbb{Z}^{2}$.
\end{lemma}

\begin{proof}
We define $\varepsilon:=$ sign$\left(  ad-bc\right)  $ and write $\kappa
c-\lambda d=\gamma q+p,$ $\gamma\in\mathbb{Z}.$ Setting $g:=\varepsilon
\kappa+\gamma b$ and $e:=\varepsilon\lambda+\gamma a,$ we get
\[
gc-ed=\varepsilon\left(  \kappa c-\lambda d\right)  +\gamma\left(
bc-ad\right)  =\varepsilon\left(  \gamma q+p\right)  +\gamma\left(
-\varepsilon q\right)  =\varepsilon p,
\]
i.e., $p=\varepsilon\left(  gc-ed\right)  $. On the other hand,
\[
\text{det}\left(
\begin{smallmatrix}
a & e\\
b & g
\end{smallmatrix}
\right)  =ag-eb=\varepsilon\left(  \kappa a-\lambda b\right)  =\varepsilon,
\]
which means that $\mathbf{n}^{\prime\prime}$ is primitive, $\left\{
\mathbf{n},\mathbf{n}^{\prime\prime}\right\}  \ $a $\mathbb{Z}$-basis of
$\mathbb{Z}^{2},$ and $\left(
\begin{smallmatrix}
a & c\\
b & d
\end{smallmatrix}
\right)  =\left(
\begin{smallmatrix}
a & e\\
b & g
\end{smallmatrix}
\right)  \left(
\begin{smallmatrix}
1 & p\\
0 & q
\end{smallmatrix}
\right)  $, i.e., $\mathbf{n}^{\prime}=p\mathbf{n}+q\mathbf{n}^{\prime\prime
},$ because%
\[
pa+qe=\varepsilon\left(  gca-eda\right)  +\varepsilon\left(  ad-bc\right)
e=c\varepsilon\left(  ga-be\right)  =c
\]
and
\[
pb+qg=\varepsilon\left(  gcb-edb\right)  +\varepsilon\left(  ad-bc\right)
g=\varepsilon d\left(  ag-be\right)  =d\ .
\]
Since gcd$\left(  p,q\right)  $ divides both $c$ and $d,$ and gcd$\left(
c,d\right)  =1$, we obtain gcd$\left(  p,q\right)  =1.$
\end{proof}

\begin{lemma}
\label{PEQU2}There is a linear map $\Phi:\mathbb{R}^{2}\longrightarrow
\mathbb{R}^{2},$ $\Phi\left(  \mathbf{x}\right)  :=\Xi\,\mathbf{x,}$ with
$\Xi\in$ \emph{GL}$_{2}(\mathbb{Z}),$ such that%
\[
\Phi\left(  \sigma\right)  =\mathbb{R}_{\geq0}\tbinom{1}{0}+\mathbb{R}_{\geq
0}\tbinom{p}{q}.
\]

\end{lemma}

\begin{proof}
It it enough to define as $\Xi:=\left(
\begin{smallmatrix}
\frac{\varepsilon\left(  d-bp\right)  }{q} & \frac{\varepsilon\left(
ap-c\right)  }{q}\\%
\genfrac{.}{.}{0pt}{}{{}}{-\varepsilon b}%
&
\genfrac{.}{.}{0pt}{}{{}}{\varepsilon a}%
\end{smallmatrix}
\right)  $.
\end{proof}

\noindent{}Henceforth, we call $\sigma$ a $(p,q)$-\textit{cone}. Denoting by
$U_{\sigma}:=$ Spec$(\mathbb{C}[\sigma^{\vee}\cap\mathbb{Z}^{2}])$ the affine
toric variety associated to $\sigma$ (by means of the monoid $\sigma^{\vee
}\cap\mathbb{Z}^{2}$, where $\sigma^{\vee}$ is the dual of $\sigma$) and by
orb$(\sigma)$ the single point being fixed under the usual action of the
algebraic torus $\mathbb{T}:=$ Hom$_{\mathbb{Z}}(\mathbb{Z}^{2},\mathbb{C}%
^{\ast})$ on $U_{\sigma},$ it is easy to see that $U_{\sigma}\cong%
\mathbb{C}^{2}$ only if $q=1.$ (In this case, $\sigma$ is said to be a
\textit{basic cone}.) On the other hand, whenever $q>1$ we have the following:

\begin{proposition}
$\emph{orb}(\sigma)\in U_{\sigma}$ is a cyclic quotient singularity. In
particular,
\[
U_{\sigma}\cong\mathbb{C}^{2}/G=\emph{Spec}(\mathbb{C}[z_{1},z_{2}]^{G}),
\]
with $G\subset$ \emph{GL}$\left(  2,\mathbb{C}\right)  $ denoting the cyclic
group $G$ of order $q$ which is generated by \emph{diag}$(\zeta_{q}^{-p}%
,\zeta_{q})$ \emph{(}$\zeta_{q}:=$ \emph{exp}$(2\pi\sqrt{-1}/q)$\emph{)} and
acts on $\mathbb{C}^{2}=$ \emph{Spec}$(\mathbb{C}[z_{1},z_{2}])$ linearly and effectively.
\end{proposition}

\begin{proof}
\noindent See \cite[\S \ 2.2, pp. 32-34]{Fulton} or \cite[Proposition 1.24,
p.30]{Oda}.$\medskip$
\end{proof}

\noindent{}In fact, $U_{\sigma}$ is the toric variety $X_{\Delta_{\sigma}}$
defined by the fan
\[
\Delta_{\sigma}:=\left\{  \sigma\ \text{together with its faces}\right\}  ,
\]
and by Proposition \ref{ISO} these two numbers $p=p_{\sigma}$ and
$\ q=q_{\sigma}$ parametrize uniquely the isomorphism class of the
germ\emph{\ }$\left(  U_{\sigma}\text{\emph{, }orb}\left(  \sigma\right)
\right)  $, up to replacement of\emph{\ }$p$ by its socius $\widehat{p}$
(which corresponds just to the interchange of the coordinates). [The
\textit{socius} $\widehat{p}$ of $p$ is defined to be the uniquely determined
integer, so that $0\leq\widehat{p}<q\emph{,}$ gcd$(\widehat{p},q)=1,$ and
$p\,\widehat{p}\equiv1$(mod $q$).]

\begin{proposition}
\label{ISO}Let $\sigma,\tau\subset\mathbb{R}^{2}$ be two $2$-dimensional,
rational, stronly convex polyhedral cones. Then the following conditions are
equivalent\emph{:} \smallskip\newline\emph{(i) \ }There is a $\mathbb{T}%
$-equivariant isomorphism $U_{\sigma}\cong U_{\tau}$ mapping \emph{orb}%
$\left(  \sigma\right)  $ onto \emph{orb}$\left(  \tau\right)  $%
.\smallskip\newline\emph{(ii)} There exists a linear map $\Phi:\mathbb{R}%
^{2}\longrightarrow\mathbb{R}^{2},$ $\Phi\left(  \mathbf{x}\right)
:=\Xi\,\mathbf{x,}$ with $\Xi\in$ \emph{GL}$_{2}(\mathbb{Z}),$ such that\emph{
}$\Phi\left(  \sigma\right)  =\tau.\smallskip$\newline\emph{(iii)} For the
numbers $p_{\sigma},$ $p_{\tau},$ $q_{\sigma},$ $q_{\tau}$ associated to
$\sigma,\tau$ \emph{(}by \emph{Lemma} \emph{\ref{pequ})} we have $q_{\tau
}=q_{\sigma}$ and either $p_{\tau}=p_{\sigma}$ or $p_{\tau}=\widehat
{p}_{\sigma}.$
\end{proposition}

\begin{proof}
For the equivalence (i)$\Leftrightarrow$(ii) see Ewald \cite[Ch. VI, Thm.
2.11, pp. 222-223]{Ewald}.\smallskip\ For proving (ii)$\Leftrightarrow$(iii)
we may w.l.o.g. consider (by virtue of Lemma \ref{PEQU2}) the cones
\[
\overline{\sigma}:=\mathbb{R}_{\geq0}\tbinom{1}{0}+\mathbb{R}_{\geq0}%
\tbinom{p_{\sigma}}{q_{\sigma}}\ \ \ \text{and\ \ \ }\overline{\tau
}:=\mathbb{R}_{\geq0}\tbinom{1}{0}+\mathbb{R}_{\geq0}\tbinom{p_{\tau}}%
{q_{\tau}}%
\]
instead of $\sigma,\tau.\smallskip$\newline(ii)$\Rightarrow$(iii): If there is
a linear map $\Phi:\mathbb{R}^{2}\longrightarrow\mathbb{R}^{2},$ $\Phi\left(
\mathbf{x}\right)  :=\Xi\,\mathbf{x,}$ with $\Xi\in$ GL$_{2}(\mathbb{Z}),$
such that\emph{ }$\Phi\left(  \overline{\sigma}\right)  =\overline{\tau},$
then either
\[
\Phi\left(  \tbinom{1}{0}\right)  =\tbinom{1}{0}\text{\ \ \ and\ \ }%
\Phi\left(  \tbinom{p_{\sigma}}{q_{\sigma}}\right)  =\tbinom{p_{\tau}}%
{q_{\tau}}%
\]
or
\[
\Phi\left(  \tbinom{1}{0}\right)  =\tbinom{p_{\tau}}{q_{\tau}}%
\text{\ \ \ and\ \ }\Phi\left(  \tbinom{p_{\sigma}}{q_{\sigma}}\right)
=\tbinom{1}{0}.
\]
Thus, either%
\[
\Xi=\left(
\begin{array}
[c]{cc}%
1 & \tfrac{p_{\tau}-p_{\sigma}}{q_{\sigma}}\smallskip\\
0 & \tfrac{q_{\tau}}{q_{\sigma}}%
\end{array}
\right)  \text{ \ \ \ \ or \ \ \ }\Xi=\allowbreak\allowbreak\left(
\begin{array}
[c]{cc}%
p_{\tau} & \tfrac{1-p_{\sigma}\ p_{\tau}}{q_{\sigma}}\smallskip\\
q_{\tau} & -\tfrac{p_{\sigma}\ q_{\tau}}{q_{\sigma}}%
\end{array}
\right)  .
\]
In the first case det$\left(  \Xi\right)  $ has to be equal to $1$, which
means that $q_{\sigma}=q_{\tau}$ and $p_{\tau}-p_{\sigma}\equiv0$ (mod
$q_{\sigma}$), i.e., $p_{\tau}=p_{\sigma}$ (because $0\leq p_{\sigma},p_{\tau
}\leq q_{\sigma}=q_{\tau}$). In the second case, det$\left(  \Xi\right)  =-1$;
hence, $q_{\sigma}=q_{\tau}$ and $1-p_{\sigma}\ p_{\tau}\equiv0$ (mod
$q_{\sigma}$), i.e., $p_{\tau}=\widehat{p}_{\sigma}$.\smallskip\ \newline%
(iii)$\Rightarrow$(ii): If $q_{\sigma}=q_{\tau}$ and $p_{\sigma}=p_{\tau}$, we
define $\Phi:=$ id$_{\mathbb{R}^{2}}$. Otherwise, $q_{\sigma}=q_{\tau}$ and
$p_{\tau}=\widehat{p}_{\sigma}$, and
\[
\Phi\left(  \mathbf{x}\right)  :=\left(
\begin{array}
[c]{cc}%
p_{\tau} & \frac{1}{q_{\sigma}}-p_{\sigma}\smallskip\\
q_{\sigma} & -p_{\sigma}%
\end{array}
\right)  \left(
\begin{array}
[c]{c}%
x_{1}\\
x_{2}%
\end{array}
\right)  ,\ \ \ \forall\mathbf{x}=\binom{x_{1}}{x_{2}}\in\mathbb{R}^{2},
\]
is an $\mathbb{R}$-vector space isomorphism with the desired property.
\end{proof}

\noindent{}To construct the minimal desingularization of $U_{\sigma}$ for a
$(p,q)$-cone
\[
\sigma=\mathbb{R}_{\geq0}\mathbf{n}+\mathbb{R}_{\geq0}\mathbf{n}^{\prime
}\subset\mathbb{R}^{2}\ \ \ \ (\text{with }q>1)
\]
we consider the negative-regular continued fraction expansion of
\[
\dfrac{q}{q-p}=\left[  \!\left[  b_{1},b_{2},\ldots,b_{s}\right]  \!\right]
:=b_{1}-\frac{1}{b_{2}-\dfrac{1}{%
\begin{array}
[c]{cc}%
\ddots & \\
& b_{s-1}-\dfrac{1}{b_{s}}%
\end{array}
}}\ \ ,
\]
and define $\mathbf{u}_{0}:=\mathbf{n},$ $\mathbf{u}_{1}:=\frac{1}%
{q}((q-p)\mathbf{n}+\mathbf{n}^{\prime}),$ and lattice points $\{\mathbf{u}%
_{j}\left\vert \,2\leq j\leq s+1\right.  \}$ by the formulae
\[
\mathbf{u}_{j+1}:=b_{j}\mathbf{u}_{j}-\mathbf{u}_{j-1},\ \ \forall
j\in\{1,\ldots,s\}.\
\]
It is easy to see that $\mathbf{u}_{s+1}=\mathbf{n}^{\prime},$ and that the
integers $b_{j}$ are $\geq2,$ for all indices $j\in\{1,\ldots,s\}.$ Next, we
subdivide $\sigma$ into $s+1$ smaller basic cones by introducing new rays
passing through the points $\mathbf{u}_{1},...,\mathbf{u}_{s}.$

\begin{theorem}
[Toric version of Hirzebruch's desingularization]\label{HIRZVER}The
refinement
\[
\widetilde{\Delta}_{\sigma}:=\left\{  \left\{  \mathbb{R}_{\geq0}%
\,\mathbf{u}_{j}+\mathbb{R}_{\geq0}\,\mathbf{u}_{j+1}\ \left\vert \ 0\leq
j\leq s\right.  \right\}  \ \ \text{\emph{together with their faces}}\right\}
\]
of $\Delta_{\sigma}:=\left\{  \sigma\ \text{\emph{together with its faces}%
}\right\}  $ consists of basic cones, is the coarsest refinement of
$\Delta_{\sigma}$ with this property, and induces the \emph{minimal}
$\mathbb{T}$-\emph{equivariant resolution} $\ X_{\widetilde{\Delta}_{\sigma}%
}\longrightarrow X_{\Delta_{\sigma}}=U_{\sigma}$ of the singular point
\emph{orb}$(\sigma).$ Moreover, the exceptional divisor is $E:=%
{\textstyle\sum\nolimits_{j=1}^{s}}
E_{j},$ having%
\[
E_{j}:=\text{ }\overline{\text{\emph{orb}}_{\widetilde{\Delta}_{\sigma}%
}(\mathbb{R}_{\geq0}\,\mathbf{u}_{j})}\ (\cong\mathbb{P}_{\mathbb{C}}%
^{1}),\ \ \forall j\in\{1,\ldots,s\},
\]
\emph{(}i.e., the closures of the $\mathbb{T}$-orbits of the new\ rays w.r.t.
$\widetilde{\Delta}_{\sigma}$\emph{)} as its components, with
self-intersection number $(E_{j})^{2}=-b_{j}.$
\end{theorem}

\begin{proof}
\noindent See Hirzebruch \cite[pp. 15-20]{Hirzebruch2} who constructs
$X_{\widetilde{\Delta}_{\sigma}}$ by resolving the unique singularity lying
over $\mathbf{0}\in\mathbb{C}^{3}$ in the normalization of the hypersurface
\[
\left\{  \left.  \left(  z_{1},z_{2},z_{3}\right)  \in\mathbb{C}%
^{3}\right\vert \,z_{1}^{q}-z_{2}z_{3}^{q-p}=0\right\}  ,
\]
and Oda \cite[pp. 24-30]{Oda} for a proof which uses only the tools of toric geometry.
\end{proof}

\section{Local indices\label{LOCALINDICES}}

\noindent Let $\sigma\subset\mathbb{R}^{2}$ be a $(p,q)$-cone. We define the
\textit{local index} $l=l_{\sigma}$ of $\sigma$ to be the positive integer%
\begin{equation}
l:=\left\{
\begin{array}
[c]{ll}%
1, & \text{if }q=1,\\
\text{min}\left\{  \left.  k\in\mathbb{N}\ \right\vert \ k\,K(E)\text{ is a
Cartier divisor}\right\}  , & \text{if }q>1,
\end{array}
\right.  \label{LOCALINDDEF}%
\end{equation}
where $K(E)$ denotes the \textit{local canonical divisor} of $X_{\widetilde
{\Delta}_{\sigma}}$ at orb$(\sigma)$ (in the sense of \cite[p. 75]{Dais})
w.r.t. the minimal resolution $X_{\widetilde{\Delta}_{\sigma}}\longrightarrow
X_{\Delta_{\sigma}}$ of orb$(\sigma)$ constructed in Theorem \ref{HIRZVER}. It
can be shown that
\begin{equation}
l=\dfrac{q}{\text{gcd}(q,p-1)}, \label{localind}%
\end{equation}
cf. \cite[Note 3.19, p. 89, and Prop. 4.4, pp. 94-95]{Dais}, and that the
self-intersection number of $K(E)$ equals%
\begin{equation}
K(E)^{2}=-\left(  \dfrac{2-\left(  p+\widehat{p}\right)  }{q}+\sum_{j=1}%
^{s}(b_{j}-2)\right)  , \label{SELFINTKE}%
\end{equation}
cf. \cite[Corollary 4.6, p. 96]{Dais}. For the proof of Theorem \ref{THM3} we
need to know under which restrictions on $p$ and $q$ we have $l\in\{1,3\}.$

\begin{lemma}
\label{lind12}If $\sigma\subset\mathbb{R}^{2}$ is a $(p,q)$-cone, then%
\begin{equation}
l=1\Longleftrightarrow\text{ }\left\{
\begin{array}
[c]{l}%
\text{\emph{either} }p=0\text{ and }q=1,\\
\text{\emph{or} }p=1\text{ and }q\geq2,
\end{array}
\right.  \label{lind1}%
\end{equation}

\end{lemma}

\begin{proof}
By (\ref{localind}), $l=1\Longleftrightarrow q=$ gcd$(q,q-p+1),$ and therefore
$q\mid p-1.$ Since $p-1<p<q,$ $p$ and $q$ satisfy conditions (\ref{lind1}%
).\hfill{}
\end{proof}

\begin{lemma}
If $\sigma\subset\mathbb{R}^{2}$ is a $(p,q)$-cone, then%
\begin{equation}
l=3\Longleftrightarrow\text{ }\left\{
\begin{array}
[c]{l}%
\text{\emph{either} }(p,q)\in A,\\
\text{\emph{or} }(p,q)\in B,
\end{array}
\right.  \label{lind3}%
\end{equation}
where%
\[
A:=\left\{  \left.  \left(  p,q\right)  \in\mathbb{N}\times\mathbb{N\ }%
\right\vert \ q=3(p-1),\ \ p\geq2,\ \ 3\nmid p\right\}  ,
\]
and%
\[
B:=\left\{  \left.  \left(  p,q\right)  \in\mathbb{N}\times\mathbb{N\ }%
\right\vert \ q=\frac{3}{2}(p-1),\ \ p\ \text{\emph{odd} }\geq5,\ \ 3\nmid
p\right\}  .
\]
Moreover, if $\left(  p,q\right)  \in A$ and $\left(  p^{\prime},q\right)  \in
B,$ then%
\[
p^{\prime}=\widehat{p}\ (=\text{\emph{the socius of} }p)\Longleftrightarrow
pp^{\prime}\equiv1(\text{\emph{mod} }q)\Longleftrightarrow q\equiv
0(\text{\emph{mod} }9).
\]

\end{lemma}

\begin{proof}
\noindent{} $l=3$ means that $q=3m,$ where $m:=$ gcd$(q,p-1).$ Write $p-1=am.$
Since $1\leq p<q,$ we have $a\in\{1,2\}.$ Since gcd$(p,q)=1,$ in the case in
which $a=1,$ we get gcd$(3m,m+1)=1\Longleftrightarrow$
gcd$(3,p)=1\Longleftrightarrow3\nmid p,$ i.e. $(p,q)\in A,$ whereas in the
case in which $a=2,$ we get gcd$(3m,2m+1)=1\Longleftrightarrow$
gcd$(3,p)=1\Longleftrightarrow3\nmid p,$ and $p$ odd $\geq5,$ i.e. $(p,q)\in
B.$ Hence, (\ref{lind3}) is true. The last assertion can be verified easily.
\end{proof}

\begin{note}
\label{NoteAB}It is worthwhile to take a closer look at the sets $A$ and $B,$
and to the corresponding negative-regular continued fraction
expansions.\smallskip\ \newline\underline{Set $A:$}\newline%
\[%
\begin{tabular}
[c]{|c|c|c|c|c|c|c|c|c|c|c|c|c|}\hline
$p$ & $2$ & $\mathit{4}$ & $5$ & $\mathit{7}$ & $8$ & $\mathit{10}$ & $11$ &
$\mathit{13}$ & $14$ & $\mathit{16}$ & $17$ & $\cdots$\\\hline\hline
$q$ & $3$ & $\mathit{9}$ & $12$ & $\mathit{18}$ & $21$ & $\mathit{27}$ & $30$
& $\mathit{36}$ & $39$ & $\mathit{45}$ & $48$ & $\cdots$\\\hline
\end{tabular}
\ \ \ \ \ \ \ \ \ \
\]
$\bullet$ First case: Whenever \fbox{$9\nmid q$} we have $\widehat{p}=p$ and%
\[
\frac{q}{q-p}=\left\{
\begin{array}
[c]{ll}%
3, & \text{if }p=2,\ q=3,\\
\  & \\
\left[  \!\left[  2,4,2\right]  \!\right]  , & \text{if }p=5,\ q=12,\\
\  & \\
\![\![2,3,\underset{\left(  \frac{q-3}{9}-2\right)  \text{\emph{-}times}%
}{\underbrace{2,...,2}},3,2]\!], & \text{if }p\geq8,\ q\geq21.
\end{array}
\right.  \medskip
\]
$\bullet$ Second case: Whenever \fbox{$9\mid q$} we have $\widehat{p}=2p-1$
and
\[
\frac{q}{q-p}=\left\{
\begin{array}
[c]{ll}%
\left[  \!\left[  2,5\right]  \!\right]  , & \text{if }p=4,\ q=9,\\
\! & \\
\![\![2,3,\underset{\left(  \frac{q}{9}-2\right)  \text{\emph{-}times}%
}{\underbrace{2,...,2}},4]\!], & \text{if }p\geq7,\ q\geq18.
\end{array}
\right.  \medskip
\]
\underline{Set $B:$}%
\[%
\begin{tabular}
[c]{|c|c|c|c|c|c|c|c|c|c|c|c|c|}\hline
$p$ & $5$ & $\mathit{7}$ & $11$ & $\mathit{13}$ & $17$ & $\mathit{19}$ & $23$
& $\mathit{25}$ & $29$ & $\mathit{31}$ & $35$ & $\cdots$\\\hline\hline
$q$ & $6$ & $\mathit{9}$ & $15$ & $\mathit{18}$ & $24$ & $\mathit{27}$ & $33$
& $\mathit{36}$ & $42$ & $\mathit{45}$ & $51$ & $\cdots$\\\hline
\end{tabular}
\ \ \ \ \ \ \ \ \ \
\]
$\bullet$ First case: Whenever \fbox{$9\nmid q$} we have $\widehat{p}=p$ and%
\[
\frac{q}{q-p}=\left\{
\begin{array}
[c]{ll}%
6, & \text{if }p=5,\ q=6,\\
\  & \\
\lbrack\![4,\underset{\left(  \frac{q-6}{9}-1\right)  \text{\emph{-}times}%
}{\underbrace{2,...,2}},4]\!], & \text{if }p\geq11,\ q\geq15.
\end{array}
\right.  \medskip
\]
$\bullet$ Second case: Whenever \fbox{$9\mid q$} we have $\widehat{p}=\frac
{1}{2}(p+1)$ and
\[
\frac{q}{q-p}=\left\{
\begin{array}
[c]{ll}%
\left[  \!\left[  5,2\right]  \!\right]  , & \text{if }p=7,\ q=9,\\
\! & \\
\![\![4,\underset{\left(  \frac{q}{9}-2\right)  \text{\emph{-}times}%
}{\underbrace{2,...,2}},3,2]\!], & \text{if }p\geq13,\ q\geq18.
\end{array}
\right.
\]
These continued fraction expansions will be useful in what follows in
\S \ref{STRATEGY}.
\end{note}

\section{Compact toric surfaces\label{COMPACTTS}}

\noindent{}Every compact toric surface is a $2$-dimensional toric variety
$X_{\Delta}$ associated to a \textit{complete} fan $\Delta$ in $\mathbb{R}%
^{2},$ i.e., a fan having $2$-dimensional cones as maximal cones and whose
support $\left\vert \Delta\right\vert $ is the entire $\mathbb{R}^{2}$ (see
\cite[Theorem 1.11, p. 16]{Oda}). Consider a complete fan $\Delta$ in
$\mathbb{R}^{2}$ and suppose that
\begin{equation}
\sigma_{i}=\mathbb{R}_{\geq0}\mathbf{n}_{i}+\mathbb{R}_{\geq0}\mathbf{n}%
_{i+1},\ \ \ i\in\{1,\ldots,\nu\}, \label{MANYCONES}%
\end{equation}
are its $2$-dimensional cones (with $\nu\geq3$ and $\mathbf{n}_{i}$ primitive
for all $i\in\{1,\ldots,\nu\}$), enumerated in such a way that $\mathbf{n}%
_{1},\ldots,\mathbf{n}_{\nu}$ go \textit{anticlockwise} around the origin
exactly once in this order (under the usual convention: $\mathbf{n}_{\nu
+1}:=\mathbf{n}_{1},$ $\mathbf{n}_{0}:=\mathbf{n}_{\nu}$). Since $\Delta$ is
simplicial, the Picard number $\rho(X_{\Delta})$ of $X_{\Delta}$ (i.e., the
rank of its Picard group Pic$(X_{\Delta})$) equals%
\begin{equation}
\rho(X_{\Delta})=\text{ }\nu-2, \label{Picardnr}%
\end{equation}
(see \cite[p. 65]{Fulton}). Now suppose that $\sigma_{i}$ is a $(p_{i},q_{i}%
)$-cone for all $i\in\{1,\ldots,\nu\}$ and introduce the notation%
\begin{equation}
I_{\Delta}:=\left\{  \left.  i\in\{1,\ldots,\nu\}\ \right\vert \ q_{i}%
>1\right\}  ,\ \ J_{\Delta}:=\left\{  \left.  i\in\{1,\ldots,\nu
\}\ \right\vert \ q_{i}=1\right\}  , \label{IJNOT}%
\end{equation}
to separate the indices corresponding to non-basic from those corresponding to
basic cones. By \cite[Theorem 1.10, p. 15]{Oda} the singular locus of
$X_{\Delta}$ equals%
\[
\text{Sing}(X_{\Delta})=\left\{  \left.  \text{orb}(\sigma_{i})\ \right\vert
\ i\in I_{\Delta}\right\}  ,
\]
and its subset
\begin{equation}
\left.  \{\text{orb}(\sigma_{i})\ \right\vert \ i\in\breve{I}_{\Delta
}\},\text{ \ with \ \ }\breve{I}_{\Delta}:=\left\{  \left.  i\in I_{\Delta
}\ \right\vert \ p_{i}=1\right\}  , \label{IHAT}%
\end{equation}
constitutes the set of the \textit{Gorenstein singularities} of $X_{\Delta}.$
For all $i\in I_{\Delta}$ write
\begin{equation}
\frac{q_{i}}{q_{i}-p_{i}}=\left[  \!\!\left[  b_{1}^{(i)},b_{2}^{(i)}%
,\ldots,b_{s_{i}}^{(i)}\right]  \!\!\right]  \label{EXPPQCF}%
\end{equation}
and, in accordance with what is already mentioned for a single $2$-dimensional
non-basic cone in \S \ref{2DIMTORSING},\ define
\[
\mathbf{u}_{1}^{(i)}:=\mathbf{n}_{i},\ \ \ \mathbf{u}_{1}^{(i)}:=\frac
{1}{q_{i}}((q_{i}-p_{i})\mathbf{n}_{i}+\mathbf{n}_{i+1}),
\]
and
\[
\mathbf{u}_{j+1}^{(i)}=b_{j}^{(i)}\mathbf{u}_{j}^{(i)}-\mathbf{u}_{j-1}%
^{(i)},\ \ \forall j\in\{1,\ldots,s_{i}\}\ \ (\text{with }\mathbf{u}_{s_{i}%
+1}^{(i)}=\mathbf{n}_{i+1}).
\]
By construction, the proper birational map $f:X_{\widetilde{\Delta}%
}\longrightarrow X_{\Delta}$ induced by the refinement%
\[
\widetilde{\Delta}:=\left\{
\begin{array}
[c]{c}%
\text{ the cones }\left\{  \left.  \sigma_{i}\ \right\vert \ i\in J_{\Delta
}\right\}  \text{ and }\\
\left\{  \left.  \mathbb{R}_{\geq0}\,\mathbf{u}_{j}^{(i)}+\mathbb{R}_{\geq
0}\,\mathbf{u}_{j+1}^{(i)}\ \right\vert \ i\in I_{\Delta},\ j\in
\{0,1,\ldots,s_{i}\}\right\}  ,\\
\text{together with their faces}%
\end{array}
\right\}  .
\]
of the fan $\Delta$ is the\textit{\ minimal desingularization} of $X_{\Delta
}.$ Defining
\[
\left\{
\begin{array}
[c]{ll}%
E_{j}^{(i)}:=\text{ }\overline{\text{orb}_{\widetilde{\Delta}}(\mathbb{R}%
_{\geq0}\,\mathbf{u}_{j}^{(i)})}, & \forall i\in I_{\Delta}\text{ \ and
\ }\forall j\in\{1,2,\ldots,s_{i}\},\\
\  & \\
\overline{C}_{i}:=\overline{\text{orb}_{\widetilde{\Delta}}(\mathbb{R}_{\geq
0}\,\mathbf{n}_{i})}, & \forall i\in\{1,2,\ldots,\nu\},
\end{array}
\right.
\]
one observes that $\overline{C}_{i}$ is the \textit{strict transform }of
$C_{i}:=\overline{\text{orb}_{\Delta}(\mathbb{R}_{\geq0}\,\mathbf{n}_{i})}$
w.r.t. $f,$
\[
E^{(i)}:=\sum_{j=1}^{s_{i}}E_{j}^{(i)}%
\]
the \textit{exceptional divisor} replacing orb$(\sigma_{i})$ via $f$ (with
$(E_{j}^{(i)})^{2}=-b_{j}^{(i)},$ $\forall i\in I_{\Delta}$ and $\forall
j\in\{1,2,\ldots,s_{i}\}$), and%
\begin{equation}
K_{X_{\widetilde{\Delta}}}-f^{\ast}K_{X_{\Delta}}=\sum_{i\in I_{\Delta
}\mathbb{r}\breve{I}_{\Delta}}K(E^{(i)}) \label{TORDISCREP}%
\end{equation}
the \textit{discrepancy divisor} w.r.t. $f.$ (By $K_{X_{\Delta}}%
,K_{X_{\widetilde{\Delta}}}$ we denote the canonical divisors of $X_{\Delta}$
and $X_{\widetilde{\Delta}},$ respectively.)

\begin{proposition}
The Picard number of $X_{\widetilde{\Delta}}$ equals
\begin{equation}
\rho(X_{\widetilde{\Delta}})=\sum_{i\in I_{\Delta}}s_{i}+\left(  \nu-2\right)
=10-K_{X_{\Delta}}^{2}-\sum_{i\in I_{\Delta}\mathbb{r}\breve{I}_{\Delta}%
}K(E^{(i)})^{2}. \label{Picardequality}%
\end{equation}

\end{proposition}

\begin{proof}
The first equality follows from (\ref{Picardnr}) and from the fact that
\[
\rho(X_{\widetilde{\Delta}})=\rho(X_{\Delta})+\sharp\{\text{exceptional prime
divisors w.r.t. }f\}.
\]
(\ref{TORDISCREP}) implies%
\[
K_{X_{\widetilde{\Delta}}}^{2}=K_{X_{\Delta}}^{2}+\sum_{i\in I_{\Delta
}\mathbb{r}\breve{I}_{\Delta}}K(E^{(i)})^{2}.
\]
Substituting this expression for $K_{X_{\widetilde{\Delta}}}^{2}$ into
Noether's formula
\[
K_{X_{\widetilde{\Delta}}}^{2}=10-\rho(X_{\widetilde{\Delta}}),
\]
we obtain the second equality of (\ref{Picardequality}).
\end{proof}

\begin{definition}
[The additional characteristic numbers $r_{i}$]\label{INTROOFRis}For every
$i\in\{1,..,\nu\}$ we introduce integers $r_{i}$ \textit{uniquely determined}
by the conditions:%
\begin{equation}
r_{i}\mathbf{n}_{i}=\left\{
\begin{array}
[c]{ll}%
\mathbf{u}_{s_{i-1}}^{(i-1)}+\mathbf{u}_{1}^{(i)}, & \text{if }i\in I_{\Delta
}^{\prime},\\
\mathbf{n}_{i-1}+\mathbf{u}_{1}^{(i)}, & \text{if\emph{\ }}i\in I_{\Delta
}^{\prime\prime},\\
\mathbf{u}_{s_{i-1}}^{(i-1)}+\mathbf{n}_{i+1}, & \text{if }i\in J_{\Delta
}^{\prime},\\
\mathbf{n}_{i-1}+\mathbf{n}_{i+1}, & \text{if }i\in J_{\Delta}^{\prime\prime},
\end{array}
\right.  \label{CONDri}%
\end{equation}
where%
\[
I_{\Delta}^{\prime}:=\left\{  \left.  i\in I_{\Delta}\ \right\vert
\ q_{i-1}>1\right\}  ,\ \ I_{\Delta}^{\prime\prime}:=\left\{  \left.  i\in
I_{\Delta}\ \right\vert \ q_{i-1}=1\right\}  ,
\]
and%
\[
J_{\Delta}^{\prime}:=\left\{  \left.  i\in J_{\Delta}\ \right\vert
\ q_{i-1}>1\right\}  ,\ \ J_{\Delta}^{\prime\prime}:=\left\{  \left.  i\in
J_{\Delta}\ \right\vert \ q_{i-1}=1\right\}  ,
\]
with $I_{\Delta},J_{\Delta}$ as in (\ref{IJNOT}).
\end{definition}

\noindent{}By \cite[Lemma 4.3]{Dais}, for $i\in\{1,\ldots,\nu\},$ $-r_{i}$ is
nothing but the self-intersection number $\overline{C}_{i}^{2}$ of the strict
transform\textit{ }$\overline{C}_{i}$ of $C_{i}$ w.r.t. $f.$ The triples
$(p_{i},q_{i},r_{i}),$ $i\in\{1,2,\ldots,\nu\},$ are used to define the
\textsc{wve}$^{2}$\textsc{c}-graph\textit{ }$\mathfrak{G}_{\Delta}.$

\begin{figure}[h]
\epsfig{file=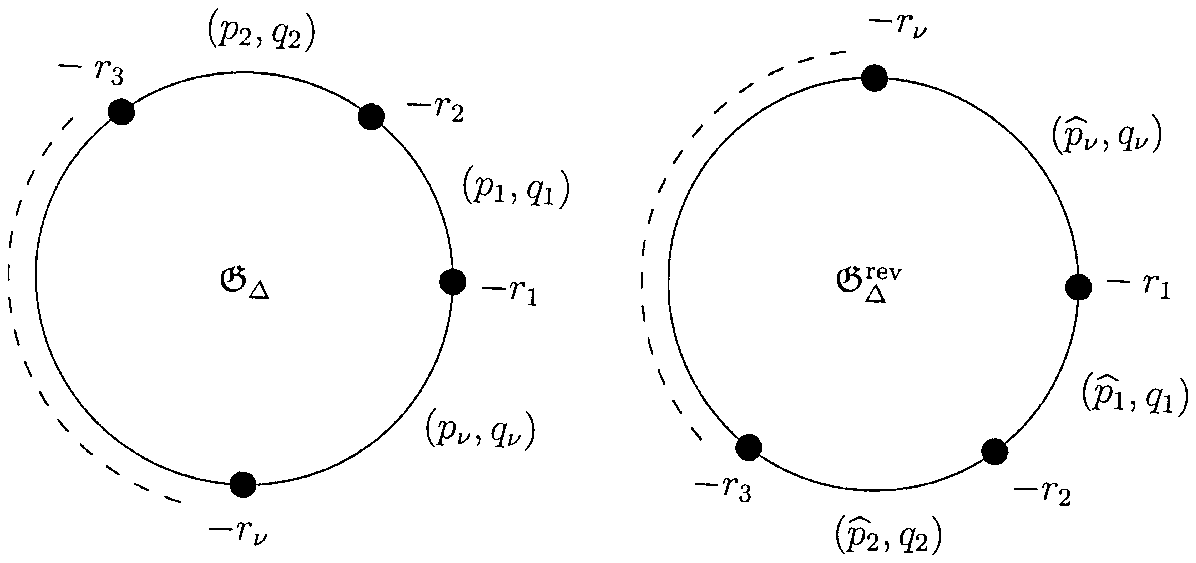, height=5cm, width=10cm}
\caption{}\label{Fig.1}%
\end{figure}

\begin{definition}
A \textit{circular graph }is a plane graph whose vertices are points on a
circle and whose edges are the corresponding arcs (on this circle, each of
which connects two consecutive vertices). We say that a circular graph
$\mathfrak{G}$ is $\mathbb{Z}$-\textit{weighted at its vertices} and
\textit{double} $\mathbb{Z}$-\textit{weighted} \textit{at its edges} (and call
it \textsc{wve}$^{2}$\textsc{c}-\textit{graph}, for short) if it is
accompanied by two maps
\[
\left\{  \text{Vertices of }\mathfrak{G}\right\}  \longmapsto\mathbb{Z}%
,\ \left\{  \text{Edges of }\mathfrak{G}\right\}  \longmapsto\mathbb{Z}^{2},
\]
assigning to each vertex an integer and to each edge a pair of integers,
respectively. To every complete fan $\Delta$ in $\mathbb{R}^{2}$ (as described
above) we associate an anticlockwise directed \textsc{wve}$^{2}$%
\textsc{c}-graph $\mathfrak{G}_{\Delta}$ with
\[
\left\{  \text{Vertices of }\mathfrak{G}_{\Delta}\right\}  =\{\mathbf{v}%
_{1},\ldots,\mathbf{v}_{\nu}\}\text{\ }\ \text{and\ \ }\left\{  \text{Edges of
}\mathfrak{G}_{\Delta}\right\}  =\{\overline{\mathbf{v}_{1}\mathbf{v}_{2}%
},\ldots,\overline{\mathbf{v}_{\nu}\mathbf{v}_{1}}\},
\]
$(\mathbf{v}_{\nu+1}:=\mathbf{v}_{1}),$ by defining its \textquotedblleft
weights\textquotedblright\ as follows:
\[
\mathbf{v}_{i}\longmapsto-r_{i},\text{ \ }\ \ \overline{\mathbf{v}%
_{i}\mathbf{v}_{i+1}}\longmapsto\left(  p_{i},q_{i}\right)  ,\ \forall
i\in\{1,\ldots,\nu\}.
\]
The \textit{reverse graph} $\mathfrak{G}_{\Delta}^{\text{rev}}$ of
$\mathfrak{G}_{\Delta}$ is the directed \textsc{wve}$^{2}$\textsc{c}-graph
which is obtained by changing the double weight $\left(  p_{i},q_{i}\right)  $
of the edge $\overline{\mathbf{v}_{i}\mathbf{v}_{i+1}}$ into $(\widehat{p}%
_{i},q_{i})$ and reversing the initial anticlockwise direction of
$\mathfrak{G}_{\Delta}$ into clockwise direction (see Figure \ref{Fig.1}).
\end{definition}

\begin{theorem}
[Classification up to isomorphism]\label{CLASSIFTHM}Let $\Delta,$
$\Delta^{\prime}$ be two complete fans in $\mathbb{R}^{2}.$ Then the following
conditions are equivalent\emph{:\smallskip\ }\newline\emph{(i)} The compact
toric surfaces $X_{\Delta}$ and $X_{\Delta^{\prime}}$ are
isomorphic.\smallskip\ \newline\emph{(ii)} Either $\mathfrak{G}_{\Delta
^{\prime}}\overset{\text{\emph{gr.}}}{\cong}\mathfrak{G}_{\Delta}$ or
$\mathfrak{G}_{\Delta^{\prime}}\overset{\text{\emph{gr.}}}{\cong}%
\mathfrak{G}_{\Delta}^{\text{\emph{rev}}}.$
\end{theorem}

\noindent Here \textquotedblleft$\,\overset{\text{gr.}}{\cong}\,$%
\textquotedblright\ indicates graph-theoretic isomorphism (i.e., a bijection
between the sets of vertices which preserves the corresponding weights). For
further details and for the proof of Theorem \ref{CLASSIFTHM} the reader is
referred to \cite[\S 5]{Dais}.

\section{Toric log Del Pezzo surfaces\label{TLDPSURFACES}}

\noindent Let $X_{\Delta}$ be a compact toric surface defined by a complete
fan $\Delta$ in $\mathbb{R}^{2}$ having (\ref{MANYCONES}) as its
$2$-dimensional cones. (Throughout this section we maintain the notation
introduced in \S \ref{COMPACTTS}.) \ It is known that $X_{\Delta}$ is a log
Del Pezzo surface if and only if the minimal generators $\mathbf{n}_{1}%
,\ldots,\mathbf{n}_{\nu}$ of the rays of $\Delta$ are vertices of a lattice
polygon $Q_{\Delta}$ (cf. \cite[Remark 6.7, p. 107]{Dais}).

\begin{definition}
A polygon $Q\subset\mathbb{R}^{2}$ is called \textit{LDP-polygon} if it
contains the origin in its interior, and its vertices are primitive elements
of $\mathbb{Z}^{2}$.
\end{definition}

\noindent In fact, there is a one-to-one correspondence
\[
\left\{
\begin{array}
[c]{c}%
\text{isomorphism classes}\\
\text{of toric log Del Pezzo }\\
\text{surfaces}%
\end{array}
\right\}  \ni\left[  X_{\Delta}\right]  \longmapsto\left[  Q_{\Delta}\right]
\in\left\{
\begin{array}
[c]{c}%
\text{lattice-equivalence }\\
\text{classes}\\
\text{of LDP-polygons}%
\end{array}
\right\}  .
\]
Indeed, if $X_{\Delta}\cong X_{\Delta^{\prime}},$ then by Theorem
\ref{CLASSIFTHM} there exists a unimodular trasformation $\Phi:\mathbb{R}%
^{2}\longrightarrow\mathbb{R}^{2}$ with $\Phi(Q_{\Delta})=Q_{\Delta^{\prime}%
}.$ The inverse of the above correspondence is given by mapping the
lattice-equivalence class $\left[  Q\right]  $ of any LDP-polygon $Q$ onto
$\left[  X_{\Delta_{Q}}\right]  ,$ where%
\[
\Delta_{Q}:=\left\{  \left.  \text{the cones }\mathbb{R}_{\geq0}%
F\,\text{\ together with their faces\ }\right\vert \,F\in\mathcal{F}%
(Q)\right\}  ,
\]
and $\mathcal{F}(Q):=\left\{  \text{facets (edges) of }Q\right\}  .$ (If \ $Q$
is an LDP-polygon,
\[
\Phi:\mathbb{R}^{2}\longrightarrow\mathbb{R}^{2},\ \ \Phi\left(
\mathbf{x}\right)  :=\Xi\,\mathbf{x,\ \forall x}\in\mathbb{R}^{2},\ \text{with
}\Xi\in\text{GL}_{2}(\mathbb{Z}),
\]
and $Q^{\prime}:=\Phi(Q),$ then $\mathfrak{G}_{\Delta_{Q^{\prime}}}%
\overset{\text{gr.}}{\cong}\mathfrak{G}_{\Delta_{Q}}$ whenever det$(\Xi)=1,$
and $\mathfrak{G}_{\Delta_{Q^{\prime}}}\overset{\text{gr.}}{\cong}%
\mathfrak{G}_{\Delta_{Q}}^{\text{rev}}$ whenever det$(\Xi)=-1.$)

Therefore, the classification of toric log Del Pezzo surfaces (up to
isomorphism) is equivalent to the classification of LDP-polygons (up to
unimodular transformations). Since the number of lattice-equivalence classes
of LDP-polygons $Q_{\Delta}$ for all those $X_{\Delta}$'s having \textit{fixed
}index $\ell$ (with $\ell$ as defined in \S \ref{INTRO}) is \textit{finite},
as it follows from results appearing in \cite{A-Br}, \cite{Bo-Bo},
\cite{Hensley} and \cite{Lag-Zieg}, it is reasonable (for any systematic
approach to the classification problem) to focus on $\ell.$ By
(\ref{LOCALINDDEF}), (\ref{localind}), (\ref{IHAT}) and (\ref{TORDISCREP}) we obtain:

\begin{lemma}
\label{indlind}The index $\ell$ of a toric log Del Pezzo surface $X_{\Delta}$
equals%
\begin{equation}
\ell=\left\{
\begin{array}
[c]{ll}%
\operatorname{lcm}\left\{  \left.  l_{i}\ \right\vert \ i\in I_{\Delta
}\right\}  \ (=\operatorname{lcm}\{\left.  l_{i}\ \right\vert \ i\in
I_{\Delta}\mathbb{r}\breve{I}_{\Delta}\}), & \text{\emph{if} \ }I_{\Delta}%
\neq\varnothing,\\
1, & \text{\emph{if} \ }I_{\Delta}=\varnothing,
\end{array}
\right.  \label{GLOBALELL}%
\end{equation}
where $l_{i}=l_{\sigma_{i}}$\emph{ }is the local index of $\sigma_{i}$\emph{
(}cf. \emph{(\ref{localind})).}
\end{lemma}

\begin{remark}
In geometric terms, $\ell=$ min$\{\left.  k\in\mathbb{N}\right\vert
\!kQ_{\Delta}^{\ast}$ is a lattice polygon$\}$, where $Q_{\Delta}^{\ast}$
denotes the \textit{polar} of the polygon $Q_{\Delta}.$ In other words, $\ell$
equals the least common multiple of the (smallest) denominators of the
(rational) coordinates of the vertices of $Q_{\Delta}^{\ast}.$ Moreover, for
$\ell\geq2,$ $\nu=\sharp\{$vertices of $Q_{\Delta}\}\leq4\ell+1$ (see
\cite[Lemma 3.1]{DN}).
\end{remark}

\begin{proposition}
For any toric log Del Pezzo surface $X_{\Delta}$ of index $\ell\geq1$ the
following inequality holds\emph{:}
\begin{equation}
\sum_{i\in I_{\Delta}}s_{i}\leq12-\sum_{i\in I_{\Delta}\mathbb{r}\breve
{I}_{\Delta}}K(E^{(i)})^{2}-\left(  1+\frac{1}{\ell}\right)  \nu.
\label{Picardinequality}%
\end{equation}

\end{proposition}

\begin{proof}
(\ref{Picardinequality}) follows from (\ref{Picardequality}) and
$K_{X_{\Delta}}^{2}\geq\frac{\nu}{\ell}$ (see \cite[proof of Lemma 3.2]{DN}).
\end{proof}

\noindent An additional necessary condition for a compact toric surface
$X_{\Delta}$ to be log Del Pezzo is dictated by the \textit{convexity} of the
necessarily existing LDP-polygon $Q_{\Delta}$:

\begin{proposition}
For any toric log Del Pezzo surface $X_{\Delta}$ of index $\ell\geq2$ we have%
\begin{equation}%
{\displaystyle\sum\limits_{i\in\breve{I}_{\Delta}}}
q_{i}\leq\left(
{\displaystyle\sum\limits_{i\in I_{\Delta}\mathbb{r}\breve{I}_{\Delta}}}
\left(  1-\frac{2}{l_{i}}\right)  q_{i}\right)  -(\nu-\sharp(I_{\Delta}))+8.
\label{Scottineq}%
\end{equation}

\end{proposition}

\begin{proof}
Since%
\[
\sharp(\text{int}(\text{conv}(\{\mathbf{n}_{i},\mathbf{n}_{i+1}\})\cap
\mathbb{Z}^{2})=\text{ gcd}(q_{i},p_{i}-1)-1,\ \ \forall i\in\{1,\ldots
,\nu\},
\]
we obtain%
\begin{equation}
\sharp(\partial Q_{\Delta}\cap\mathbb{Z}^{2})=\nu+\sum_{i=1}^{\nu}%
\,\sharp(\text{int}(\text{conv}(\{\mathbf{n}_{i},\mathbf{n}_{i+1}%
\})\cap\mathbb{Z}^{2})=\sum_{i=1}^{\nu}\text{gcd}(q_{i},p_{i}-1).
\label{boundary}%
\end{equation}
($\partial,$ int, and conv are used as abbreviations for boundary, interior,
and convex hull, respectively.) Furthermore, since
\[
\text{area}(Q_{\Delta})=\sum_{i=1}^{\nu}\text{area}(\text{conv}(\{\mathbf{0}%
,\mathbf{n}_{i},\mathbf{n}_{i+1}\}))=\frac{1}{2}\left(  \sum\limits_{i=1}%
^{\nu}q_{i}\right)  ,
\]
using Pick's formula (cf. \cite[p. 113]{Fulton}, \cite[p. 101]{Oda}):%
\[
\sharp(Q_{\Delta}\cap\mathbb{Z}^{2})=\text{area}(Q_{\Delta})+\frac{1}{2}%
\sharp(\partial Q_{\Delta}\cap\mathbb{Z}^{2})+1,
\]
we get%
\begin{equation}
\sharp(\text{int}(Q_{\Delta})\cap\mathbb{Z}^{2})=\frac{1}{2}\left(  \sum
_{i=1}^{\nu}(q_{i}-\text{gcd}(q_{i},p_{i}-1))\right)  +1. \label{interior}%
\end{equation}
Finally, since $\ell\geq2,$ Scott's inequality \cite{Scott} can be written as%
\begin{equation}
\sharp(\partial Q_{\Delta}\cap\mathbb{Z}^{2})<2\,\sharp(\text{int}(Q_{\Delta
})\cap\mathbb{Z}^{2})+7. \label{Scottoriginal}%
\end{equation}
By (\ref{boundary}), (\ref{interior}), (\ref{Scottoriginal}) and
(\ref{localind}) we infer that%
\[
\sum_{i=1}^{\nu}\left(  \frac{2}{l_{i}}-1\right)  q_{i}\leq8,
\]
which can be rewritten (by keeping the involved $q_{i}$'s with non-negative
coefficients) in the form (\ref{Scottineq}).
\end{proof}

\section{Compact toric surfaces with Picard number 1\label{CTSWITHPIC1}}

{}\noindent By virtue of (\ref{Picardnr}) the compact toric surfaces with
Picard number $1$ are defined by complete fans $\Delta$ in $\mathbb{R}^{2}$
\noindent{}\noindent{}with exactly three $2$-dimensional cones. Let $\Delta$
be a complete fan of this kind and%
\begin{equation}
\sigma_{1}=\mathbb{R}_{\geq0}\mathbf{n}_{1}+\mathbb{R}_{\geq0}\mathbf{n}%
_{2},\ \ \sigma_{2}=\mathbb{R}_{\geq0}\mathbf{n}_{2}+\mathbb{R}_{\geq
0}\mathbf{n}_{3},\ \ \sigma_{3}=\mathbb{R}_{\geq0}\mathbf{n}_{3}%
+\mathbb{R}_{\geq0}\mathbf{n}_{1}, \label{3CONES}%
\end{equation}
be its $2$-dimensional cones, with $\mathbf{n}_{i}$ primitive and $\sigma_{i}$
a $(p_{i},q_{i})$-cone for $i\in\{1,2,3\}.$

\begin{lemma}
\label{WEIGHTEDLEMMA}$X_{\Delta}$ is isomorphic to the quotient space
$\mathbb{P}_{\mathbb{C}}^{2}(q_{1},q_{2},q_{3})/H_{\Delta},$ where $H_{\Delta
}$ is a finite abelian group of order \emph{gcd}$(q_{1},q_{2},q_{3}).$
\end{lemma}

\begin{proof}
Since $q_{i}=\left\vert \det(\mathbf{n}_{i},\mathbf{n}_{i+1})\right\vert $ for
$i\in\{1,2,3\},$ using Cramer's rule we obtain
\[
q_{1}\mathbf{n}_{3}+q_{2}\mathbf{n}_{1}+q_{3}\mathbf{n}_{2}=\mathbf{0}.
\]
By \cite[Proposition 4.7, p. 224]{Conrads} we have $X_{\Delta}\cong%
\mathbb{P}_{\mathbb{C}}^{2}(q_{1},q_{2},q_{3})/H_{\Delta},$ where $H_{\Delta}$
is a group isomorphic to $\mathbb{Z}^{2}/(\oplus_{i=1}^{3}\mathbb{Z}%
\mathbf{n}_{i}).$ By
\[
\left\vert H_{\Delta}\right\vert =\sharp\left(  \{\text{fundamental
perallelepiped of }\oplus_{i=1}^{3}\mathbb{Z}\mathbf{n}_{i}\}\cap
\mathbb{Z}^{2}\right)  =\det(\oplus_{i=1}^{3}\mathbb{Z}\mathbf{n}_{i}),
\]
and the fact that $\det(\oplus_{i=1}^{3}\mathbb{Z}\mathbf{n}_{i})=$
gcd$(q_{1},q_{2},q_{3}),$ the assertion is true.
\end{proof}

\begin{figure}[h]
\epsfig{file=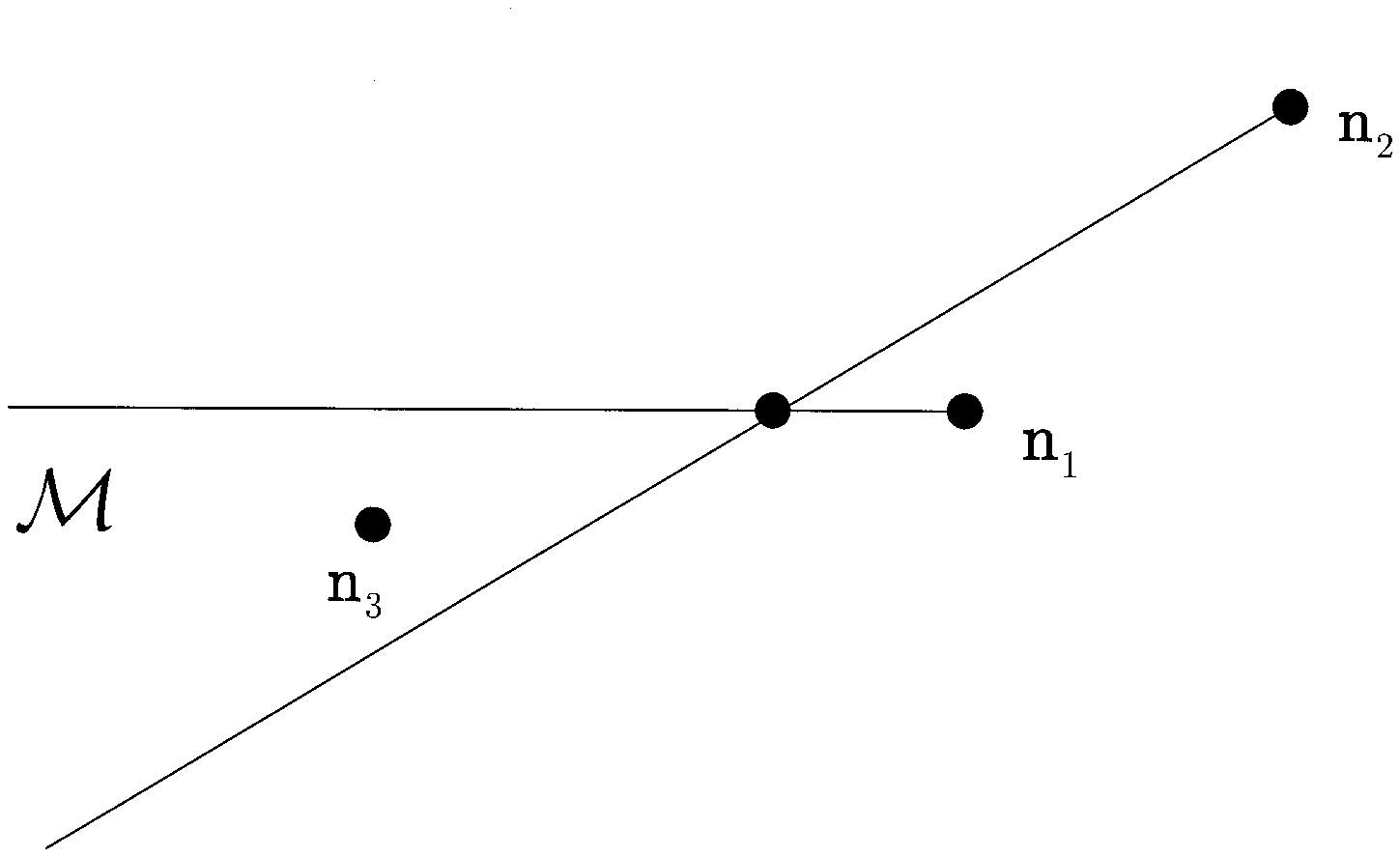, height=4.5cm, width=8cm}
\caption{}\label{Fig.2}%
\end{figure}

\noindent{}Since we are interested in describing $X_{\Delta}$ \textit{up to
isomorphism} (cf. Lemma \ref{PEQU2} and Theorem \ref{CLASSIFTHM}) we may
henceforth assume, without loss of generality, that $\mathbf{n}_{1}=\binom
{1}{0}$ and $\mathbf{n}_{2}=\binom{p_{1}}{q_{1}}.$ As all cones of $\Delta$
are strongly convex, $\mathbf{n}_{3}$ belongs (as shown in Figure \ref{Fig.2})
necessarily to the set%
\[
\mathcal{M}:=\left\{  \tbinom{x}{y}\in\mathbb{Z}^{2}\ \left\vert
\ \tfrac{q_{1}}{p_{1}}x<y<0\right.  \right\}  .
\]

\begin{lemma}
\label{MAINLEMMA1}We have
\begin{equation}
\mathbf{n}_{3}=\tbinom{-(q_{2}+p_{1}q_{3})/q_{1}}{-q_{3}}, \label{N3}%
\end{equation}
and therefore $q_{1}\mid q_{2}+p_{1}q_{3}$ and \emph{gcd}$((q_{2}+p_{1}%
q_{3})/q_{1},q_{3})=1.$ Moreover,
\begin{equation}
q_{1}q_{2}\mid\widehat{p}_{1}q_{2}+p_{2}q_{1}+q_{3}, \label{DIVCON1}%
\end{equation}
and
\begin{equation}
q_{1}q_{3}\mid p_{1}q_{3}+\widehat{p}_{3}q_{1}+q_{2}. \label{DIVCON2}%
\end{equation}

\end{lemma}

\begin{proof}
We use Lemma \ref{PEQU2}. Since $\sigma_{2}$ is a $(p_{2},q_{2})$-cone and
$\sigma_{3}$ is a $(p_{3},q_{3})$-cone, \noindent{}setting $\mathbf{n}%
_{3}=\binom{x}{y},$ we have $\ $%
\begin{equation}
\left.
\begin{array}
[c]{r}%
\left\vert \det\left(
\begin{smallmatrix}
x & p_{1}\\
y & q_{1}%
\end{smallmatrix}
\right)  \right\vert =q_{2},\\
\\
\binom{x}{y}\in\mathcal{M}%
\end{array}
\right\}  \Longrightarrow q_{1}x-p_{1}y=-q_{2}. \label{DETCON1}%
\end{equation}
on the one hand, and%
\[
\left.
\begin{array}
[c]{r}%
\left\vert \det\left(
\begin{smallmatrix}
x & 1\\
y & 0
\end{smallmatrix}
\right)  \right\vert =q_{3},\\
\\
\binom{x}{y}\in\mathcal{M}%
\end{array}
\right\}  \Longrightarrow y=-q_{3},
\]
on the other. Hence, (\ref{DETCON1}) gives $x=-\frac{1}{q_{1}}(q_{2}%
+p_{1}q_{3}).$ Moreover, by the definition of $\widehat{p}_{1}$ there exists
an integer $\lambda$ such that%
\[
\widehat{p}_{1}p_{1}-\lambda q_{1}=1.
\]
This means that
\[
\widehat{p}_{1}(-\tfrac{1}{q_{1}}(q_{2}+p_{1}q_{3}))-\lambda(-q_{3})\equiv
p_{2}(\text{mod }q_{2}),
\]
i.e., there is a $\mu\in\mathbb{Z}$ with $\mu q_{2}=p_{2}+\frac{1}{q_{1}%
}\left(  \widehat{p}_{1}(q_{2}+p_{1}q_{3})-\lambda q_{3}q_{1}\right)  .$
Consequently,
\[
\mu q_{1}q_{2}=\widehat{p}_{1}q_{2}+q_{3}(\widehat{p}_{1}p_{1}-\lambda
q_{1})+p_{2}q_{1}=\widehat{p}_{1}q_{2}+p_{2}q_{1}+q_{3},
\]
$\mu\in\mathbb{N},$ and the divisibility condition (\ref{DIVCON1}) is true.
Next, by Lemma \ref{PEQU2} there is a matrix $\left(
\begin{smallmatrix}
\mathfrak{a} & \mathfrak{b}\\
\mathfrak{c} & \mathfrak{d}%
\end{smallmatrix}
\right)  $ $\in$ GL$_{2}(\mathbb{Z})$ such that $\left(
\begin{smallmatrix}
\mathfrak{a} & \mathfrak{b}\\
\mathfrak{c} & \mathfrak{d}%
\end{smallmatrix}
\right)  \tbinom{x}{y}=\tbinom{1}{0}$ and $\ \left(
\begin{smallmatrix}
\mathfrak{a} & \mathfrak{b}\\
\mathfrak{c} & \mathfrak{d}%
\end{smallmatrix}
\right)  \tbinom{1}{0}=\tbinom{p_{3}}{q_{3}},$ i.e., $\mathfrak{a}=p_{3},$
$\mathfrak{c}=q_{3},$ and\smallskip%
\[
\left\{
\begin{array}
[c]{l}%
q_{3}x+\mathfrak{d}y=q_{3}x-\mathfrak{d}q_{3}=0\Longrightarrow\mathfrak{d}%
=x,\\
\left.  \left.
\begin{array}
[c]{c}%
p_{3}x+\mathfrak{b}y=p_{3}x-\mathfrak{b}q_{3}=1\\
x<0
\end{array}
\right\}  \right.  \Longrightarrow x=\widehat{p}_{3}-\kappa q_{3},\text{ for
some }\kappa\in\mathbb{N}.
\end{array}
\right.
\]
By (\ref{DETCON1}),
\[
q_{1}x-p_{1}y=q_{1}\left(  \widehat{p}_{3}-\kappa q_{3}\right)  +p_{1}%
q_{3}=-q_{2}\Longrightarrow\kappa q_{1}q_{3}=p_{1}q_{3}+\widehat{p}_{3}%
q_{1}+q_{2},
\]
leading to the divisibility condition (\ref{DIVCON2}).
\end{proof}

The converse is also true.

\begin{lemma}
\label{MAINLEMMA2}Given a triple of pairs $\left\{  \left.  (p_{i}%
,q_{i})\right\vert 1\leq i\leq3\right\}  $ of non-negative integers with
$p_{i}<q_{i}$ and $\gcd(p_{i},q_{i})=1$ for $i\in\{1,2,3\}$, and such that
\[
q_{1}q_{2}\mid\widehat{p}_{1}q_{2}+p_{2}q_{1}+q_{3}\ \ \ \text{and}%
\ \ \ \ q_{1}q_{3}\mid p_{1}q_{3}+\widehat{p}_{3}q_{1}+q_{2},
\]
the $2$-dimensional cones%
\[
\sigma_{1}=\mathbb{R}_{\geq0}\tbinom{1}{0}+\mathbb{R}_{\geq0}\tbinom{p_{1}%
}{q_{1}},\ \ \sigma_{2}=\mathbb{R}_{\geq0}\tbinom{p_{1}}{q_{1}}+\mathbb{R}%
_{\geq0}\tbinom{-(q_{2}+p_{1}q_{3})/q_{1}}{-q_{3}},
\]%
\[
\text{and \ }\sigma_{3}=\mathbb{R}_{\geq0}\tbinom{-(q_{2}+p_{1}q_{3})/q_{1}%
}{-q_{3}}+\mathbb{R}_{\geq0}\tbinom{1}{0},
\]
\emph{(}written by means of their minimal generators\emph{)} compose, together
with their faces, a complete fan in $\mathbb{R}^{2}$ and $\sigma_{i}$ is a
$(p_{i},q_{i})$-cone, for $i\in\{1,2,3\}.$
\end{lemma}

\begin{proof}
Obviously, $\sigma_{1}$ is a $(p_{1},q_{1})$-cone and%
\[
\det\left(
\begin{array}
[c]{cc}%
p_{1} & -(q_{2}+p_{1}q_{3})/q_{1}\\
q_{1} & -q_{3}%
\end{array}
\right)  =q_{2},\ \ \ \det\left(
\begin{array}
[c]{cc}%
-(q_{2}+p_{1}q_{3})/q_{1} & 1\\
-q_{3} & 0
\end{array}
\right)  =q_{3}.
\]
Furthermore,
\[
q_{1}q_{3}\mid p_{1}q_{3}+\widehat{p}_{3}q_{1}+q_{2}\Longrightarrow q_{1}\mid
q_{2}+p_{1}q_{3}\Longrightarrow\tbinom{-(q_{2}+p_{1}q_{3})/q_{1}}{-q_{3}}%
\in\mathbb{Z}^{2},
\]
and setting $\delta:=$ gcd$(q_{2}+p_{1}q_{3},q_{1}q_{3})$ we obtain%
\[
\delta\mid p_{1}q_{3}+\widehat{p}_{3}q_{1}+q_{2}\Longrightarrow\delta
\mid\widehat{p}_{3}q_{1}\Longrightarrow\delta\mid\widehat{p}_{3}p_{3}q_{1}.
\]
Since there exists an integer $\gamma$ with $\widehat{p}_{3}p_{3}-\gamma
q_{3}=1,$ we have
\[
\left.
\begin{array}
[c]{r}%
\delta\mid\left(  \gamma q_{3}+1\right)  q_{1}\\
\delta\mid q_{3}q_{1}\Longrightarrow\delta\mid\gamma q_{3}q_{1}%
\end{array}
\right\}  \Longrightarrow\delta\mid q_{1}.
\]
This divisibility condition is equivalent to: gcd$(\frac{1}{q_{1}}(q_{2}%
+p_{1}q_{3}),q_{3})=1,$ and therefore $\binom{-(q_{2}+p_{1}q_{3})/q_{1}%
}{-q_{3}}$ is primitive. On the other hand,
\[
q_{1}q_{2}\mid\widehat{p}_{1}q_{2}+p_{2}q_{1}+q_{3}\Longrightarrow\exists
\,\mu\in\mathbb{N}:\mu q_{1}q_{2}=\widehat{p}_{1}q_{2}+p_{2}q_{1}%
+q_{3}.\smallskip
\]
Since there exists an integer $\lambda$ with $\widehat{p}_{1}p_{1}-\lambda
q_{1}=1,$ and%
\[
\mu q_{1}q_{2}=\widehat{p}_{1}q_{2}+q_{3}(\widehat{p}_{1}p_{1}-\lambda
q_{1})+p_{2}q_{1}\Longrightarrow\widehat{p}_{1}(-\tfrac{1}{q_{1}}(q_{2}%
+p_{1}q_{3}))-\lambda(-q_{3})\equiv p_{2}(\text{mod }q_{2}),
\]
$\sigma_{2}$ is a $(p_{2},q_{2})$-cone. Finally,
\[
q_{1}q_{3}\mid p_{1}q_{3}+\widehat{p}_{3}q_{1}+q_{2}\Longrightarrow
\exists\,\kappa\in\mathbb{N}:\kappa q_{1}q_{3}=p_{1}q_{3}+\widehat{p}_{3}%
q_{1}+q_{2},\text{ i.e.,}%
\]%
\[
q_{1}\left(  \widehat{p}_{3}-\kappa q_{3}\right)  +p_{1}q_{3}=-q_{2}%
=q_{1}(-\tfrac{1}{q_{1}}(q_{2}+p_{1}q_{3}))+p_{1}q_{3}\Longrightarrow
-\tfrac{1}{q_{1}}(q_{2}+p_{1}q_{3})=\widehat{p}_{3}-\kappa q_{3},
\]
giving%
\[
\left(
\begin{array}
[c]{cc}%
p_{3} & \tfrac{1}{q_{3}}\left(  p_{3}\widehat{p}_{3}-1\right)  -\kappa p_{3}\\
q_{3} & \widehat{p}_{3}-\kappa q_{3}%
\end{array}
\right)  \binom{-\tfrac{1}{q_{1}}(q_{2}+p_{1}q_{3})}{-q_{3}}=\left(
\begin{array}
[c]{c}%
1\\
0
\end{array}
\right)  ,
\]
and%
\[
\left(
\begin{array}
[c]{cc}%
p_{3} & \tfrac{1}{q_{3}}\left(  p_{3}\widehat{p}_{3}-1\right)  -\kappa p_{3}\\
q_{3} & \widehat{p}_{3}-\kappa q_{3}%
\end{array}
\right)  \left(
\begin{array}
[c]{c}%
1\\
0
\end{array}
\right)  =\left(
\begin{array}
[c]{c}%
p_{3}\\
q_{3}%
\end{array}
\right)  .
\]
Hence, as it is explained in the proof of Proposition \ref{ISO}, the cone
$\sigma_{3}$ has to be a $(p_{3},q_{3})$-cone.
\end{proof}

\begin{lemma}
\label{PIC1}Every compact toric surface $X_{\Delta}$ having Picard number
$\rho(X_{\Delta})=1$ is a log Del Pezzo surface.
\end{lemma}

\begin{proof}
If $X_{\Delta}$ is a compact toric surface with $\rho(X_{\Delta})=1,$ then the
minimal generators $\mathbf{n}_{1},\mathbf{n}_{2},\mathbf{n}_{3}$ of the tree
cones (\ref{3CONES}) of $\Delta$ have to be in general position because the
cones are strongly convex. Hence, conv$(\{\mathbf{n}_{1},\mathbf{n}%
_{2},\mathbf{n}_{3}\})$ has to be an LDP-triangle.
\end{proof}

\begin{note}
Compact toric surfaces $X_{\Delta}$ having Picard number $\rho(X_{\Delta}%
)\geq2$ are not always log Del Pezzo surfaces. For instance, the smooth
compact surfaces $X_{\Delta}$ with $\rho(X_{\Delta})=2$ are the Hirzebruch
surfaces $\mathbb{F}_{\kappa},$ $\kappa\geq0$ (cf. \cite[Corollary 1.29, p.
45]{Oda}); among them, only $\mathbb{F}_{0}\cong\mathbb{P}_{\mathbb{C}}%
^{1}\times\mathbb{P}_{\mathbb{C}}^{1}$ and $\mathbb{F}_{1}$ (i.e., a
$\mathbb{P}_{\mathbb{C}}^{2}$ blown up at one point) are Del Pezzo surfaces
(see \cite[Proposition 2.21, p. 88]{Oda} or \cite[Theorem V.8.2, p.
192]{Ewald}). The geometric reason for that is actually very simple: Since
$\mathbb{F}_{\kappa}$ can be viewed as the toric surface associated to the fan
having $\binom{0}{1},\binom{1}{0},\binom{0}{-1}$ and $\binom{-1}{\kappa}$ as
minimal generators of its rays, setting $\mathbf{T}:=\text{conv}\left(  \left\{  \tbinom{1}{0},\tbinom{0}{-1}%
,\tbinom{-1}{\kappa}\right\}  \right)$
we see that $\tbinom{0}{1}\in\partial\,\mathbf{T}$ for $\kappa=2,$ and
$\tbinom{0}{1}\in$ int$(\mathbf{T)}$ for $\kappa\geq3.$
\end{note}

\section{Classification strategy for $\rho\left(  X_{\Delta}\right)  =1$ and
$\ell=3$\label{STRATEGY}}

\begin{definition}
We call a triple of pairs%
\begin{equation}
\left\{  \left.  (p_{i},q_{i})\in\mathbb{Z}^{2}\right\vert 1\leq
i\leq3\right\}  ,\ 0\leq p_{i}<q_{i},\text{ with }\gcd(p_{i},q_{i})=1,\forall
i\in\{1,2,3\}, \label{TRIPLES}%
\end{equation}
\textit{admissible} whenever it satisfies both divisibility conditions%
\begin{equation}
\fbox{$\ q_{1}q_{2}\mid\widehat{p}_{1}q_{2}+p_{2}q_{1}+q_{3}\ $}
\label{Condition1}%
\end{equation}
and%
\begin{equation}
\fbox{$\ q_{1}q_{3}\mid p_{1}q_{3}+\widehat{p}_{3}q_{1}+q_{2}\ $}
\label{Condition2}%
\end{equation}

\end{definition}

\noindent{}To classify all toric log Del Pezzo surfaces $X_{\Delta}$ having
Picard number $1$ and index $\ell=3$ \textit{up to isomorphism} it suffices
(by Lemmas \ref{MAINLEMMA1}, \ref{MAINLEMMA2}, and \ref{PIC1}, and Theorem
\ref{CLASSIFTHM}) to determine all admissible triples of pairs, and
consequently the fans $\Delta$ having%
\[
\sigma_{1}=\mathbb{R}_{\geq0}\mathbf{n}_{1}+\mathbb{R}_{\geq0}\mathbf{n}%
_{2},\ \ \sigma_{2}=\mathbb{R}_{\geq0}\mathbf{n}_{2}+\mathbb{R}_{\geq
0}\mathbf{n}_{3},\ \ \sigma_{3}=\mathbb{R}_{\geq0}\mathbf{n}_{3}%
+\mathbb{R}_{\geq0}\mathbf{n}_{1},
\]
as $2$-dimensional cones, with $\mathbf{n}_{1}=\tbinom{1}{0},\ \mathbf{n}%
_{2}=\tbinom{p_{1}}{q_{1}},\ \mathbf{n}_{3}=\tbinom{-(q_{2}+p_{1}q_{3})/q_{1}%
}{-q_{3}}$ as minimal generators, and $Q_{\Delta}=$ conv$\left(  \left\{
\mathbf{n}_{1},\mathbf{n}_{2},\mathbf{n}_{3}\right\}  \right)  $ as their
LDP-polygons, so that%
\begin{equation}
l_{i}=l_{\sigma_{i}}\in\{1,3\},\forall i\in\{1,2,3\},\text{ and \ }%
l_{k}=3\text{ \ for at least one \ }k\in\{1,2,3\}, \label{L13}%
\end{equation}
(see (\ref{GLOBALELL})). From now on we may assume w.l.o.g \ that $l_{1}=3.$
We also keep in mind the two auxiliary conditions
\begin{equation}%
\begin{array}
[c]{c}%
\begin{array}
[c]{c}%
\frame{$%
\begin{array}
[c]{ccc}
&  & \\
& s_{1}+s_{2}+s_{3}\,\leq-%
{\displaystyle\sum\limits_{i\in I_{\Delta}\mathbb{r}\breve{I}_{\Delta}}}
K(E^{(i)})^{2}+8 & \\
&  &
\end{array}
$}%
\end{array}
\end{array}
\label{Condition3}%
\end{equation}
(where, for our convenience, we set $s_{i}:=0$ for $i\in J_{\Delta},$ cf.
(\ref{IJNOT})), and%
\begin{equation}%
\begin{array}
[c]{c}%
\frame{$%
\begin{array}
[c]{ccc}
&  & \\
&
{\displaystyle\sum\limits_{i\in\breve{I}_{\Delta}}}
q_{i}\leq\frac{1}{3}%
{\displaystyle\sum\limits_{i\in I_{\Delta}\mathbb{r}\breve{I}_{\Delta}}}
q_{i}+\sharp(I_{\Delta})+5 & \\
&  &
\end{array}
$}%
\end{array}
\label{Condition4}%
\end{equation}
(following from (\ref{Picardinequality}) and (\ref{Scottineq}), respectively,
for $\nu=\ell=3$) which have to be satisfied because of Lemma \ref{PIC1}. By
assumption, each pair $(p_{i},q_{i})$ (belonging to a triple (\ref{TRIPLES})
which will be considered as \textquotedblleft candidate\textquotedblright\ for
being admissible) is necessarily of a specific \textit{type}. All possible
types are determined by conditions (\ref{L13}), (\ref{lind1}) and
(\ref{lind3}), and are listed in Table \ref{Table1}. (Since $l_{1}=3,$
$(p_{1},q_{1})$ can be of type \thinspace$\mathbf{1},\mathbf{2},\mathbf{3}%
,\mathbf{4}$ or $\mathbf{5}$.)
\newpage
\begin{table}[h]
\setlength\extrarowheight{2pt}
{\small
\[
\ \
\begin{array}
[c]{c}%
\ \
\begin{tabular}
[c]{|c|c|c|c|c|c|}\hline
\textbf{Types} & $p_{i}$ & $\widehat{p}_{i}$ & $q_{i}$ & $s_{i}$ &
$-K(E^{(i)})^{2}$\\\hline\hline
$\mathbf{1}$ & $2$ & $2$ & $3$ & $1$ & $%
\begin{array}
[c]{c}%
\frac{1}{3}\smallskip
\end{array}
$\\\hline
$\mathbf{2}$ & $3\xi_{i}+2$ & $p_{i}$ & $9\xi_{i}+3$ & $\xi_{i}+2$ & $%
\begin{array}
[c]{c}%
\frac{4}{3}\smallskip
\end{array}
$\\\hline
$\mathbf{3}$ & $3\xi_{i}+1$ & $2p_{i}-1\ (=6\xi_{i}+1)$ & $9\xi_{i}$ &
$\xi_{i}+1$ & $\begin{array}
[c]{c}%
2 \smallskip
\end{array}$\\\hline
$\mathbf{4}$ & $6\xi_{i}+5$ & $p_{i}$ & $9\xi_{i}+6$ & $\xi_{i}+1$ & $%
\begin{array}
[c]{c}%
\frac{8}{3}\smallskip
\end{array}
$\\\hline
$\mathbf{5}$ & $6\xi_{i}+1$ & $%
\begin{array}
[c]{c}%
\frac{1}{2}(p_{i}+1)\ (=3\xi_{i}+1)\smallskip
\end{array}
$ & $9\xi_{i}$ & $\xi_{i}+1$ & $2$\\\hline
$\mathbf{6}$ & $1$ & $1$ & $\geq2$ & $q_{i}-1$ & $\begin{array}
[c]{c}%
0 \smallskip
\end{array}$\\\hline
$\mathbf{7}$ & $0$ & $0$ & $1$ & $0$ & -----\\\hline
\end{tabular}
\end{array}
\ \ \ \ \ \ \ \ \ \
\]
} \smallskip\caption{{}}%
\setlength\extrarowheight{-2pt}
\label{Table1}%
\end{table}
\vspace{-0.2cm}

\noindent Here, $\xi_{i}$ denotes an integer which is \textit{positive} for
types $\mathbf{2}$, $\mathbf{3}$ and $\mathbf{5},$ and \textit{non-negative}
for type $\mathbf{4}$. (In particular, the entries of the last two columns are
computed by the continued fraction expansions mentioned in Note \ref{NoteAB}
and by the formula (\ref{SELFINTKE}).) Although the pairs $(p_{i},q_{i})$ of
type $\mathbf{2}$ (resp., of type $\mathbf{3},\mathbf{4},\mathbf{5}$ or
$\mathbf{6}$) are \textit{infinitely many}, conditions (\ref{Condition1}),
(\ref{Condition2}), (\ref{Condition3}) and (\ref{Condition4}) force the
testable triples of pairs (\ref{TRIPLES}) to be admissible only in
\textit{finitely many} cases.

\begin{note}
If orb$(\sigma_{2})$ is a non-Gorenstein singularity, then (\ref{Condition1})
implies%
\begin{equation}
\fbox{$\ \left[  \widehat{p}_{1}q_{2}+p_{2}q_{1}+q_{3}\right]  _{9}=0\ $}
\label{MOD9I}%
\end{equation}
(where $\left[  t\right]  _{9}$ denotes the remainder in the division of a
$t\in\mathbb{Z}$ by $9$) because $3\mid q_{1}$ and $3\mid q_{2}.$ Analogously,
if orb$(\sigma_{3})$ is a non-Gorenstein singularity, then (\ref{Condition2})
implies
\begin{equation}
\fbox{$\ \left[  p_{1}q_{3}+\widehat{p}_{3}q_{1}+q_{2}\right]  _{9}=0\ $}
\label{MOD9II}%
\end{equation}
These weaker, necessary conditions (\ref{MOD9I}) and (\ref{MOD9II}) turn out
to be very useful in proving that several triples of pairs (\ref{TRIPLES}) are
not admissible.
\end{note}

\noindent{}The proof of Theorem \ref{THM3} will follow in four
steps:\smallskip

\noindent$\blacktriangleright$ \underline{\textbf{Step 1}}: We determine which
of the triples of pairs (\ref{TRIPLES}) corresponding to the $125$ $(=5^{3})$
possible type combinations $(\alpha_{1},\alpha_{2},\alpha_{3}),$ with
$\alpha_{1},\alpha_{2},\alpha_{3}\in\{\mathbf{1},\mathbf{2},\mathbf{3}%
,\mathbf{4},\mathbf{5}\},$ are admissible, i.e., those $X_{\Delta}$'s with
exactly \textit{three} non-Gorenstein singularities.\smallskip\

\noindent{}$\blacktriangleright$ \underline{\textbf{Step 2}}: We determine
which of the triples of pairs (\ref{TRIPLES}) corresponding to the $100$
$(=2\cdot(5^{2}\cdot2))$ type combinations $(\alpha_{1},\alpha_{2},\alpha
_{3}),$ with $\alpha_{1}\in\{\mathbf{1},\mathbf{2},\mathbf{3},\mathbf{4}%
,\mathbf{5}\}$ and%
\[
(\alpha_{2},\alpha_{3})\in\left(  \{\mathbf{1},\mathbf{2},\mathbf{3}%
,\mathbf{4},\mathbf{5}\}\times\{\mathbf{6},\mathbf{7}\}\right)  \cup\left(
\{\mathbf{6},\mathbf{7}\}\times\{\mathbf{1},\mathbf{2},\mathbf{3}%
,\mathbf{4},\mathbf{5}\}\right)  ,
\]
are admissible, i.e., those $X_{\Delta}$'s with exactly \textit{two}
non-Gorenstein singularities. \smallskip\

\noindent{}$\blacktriangleright$ \underline{\textbf{Step 3}}: We do the same
for the triples of pairs (\ref{TRIPLES}) corresponding to the $20$ type
combinations $(\alpha_{1},\alpha_{2},\alpha_{3}),$ with $\alpha_{1}%
\in\{\mathbf{1},\mathbf{2},\mathbf{3},\mathbf{4},\mathbf{5}\}$ and $\alpha
_{2},\alpha_{3}\in\{\mathbf{6},\mathbf{7}\},$ i.e., for those $X_{\Delta}$'s
with exactly \textit{one} non-Gorenstein singularity. \smallskip\

\noindent{}$\blacktriangleright$ \underline{\textbf{Step 4}}: We find out the
\textsc{wve}$^{2}$\textsc{c}-graphs $\mathfrak{G}_{\Delta}$ for those
$X_{\Delta}$'s determined in steps $1$-$3,$ and then, using Theorem
\ref{CLASSIFTHM}, we pick out a suitable, minimal set of representatives of
$X_{\Delta}$'s all of whose members are pairwise non-isomorphic. Finally, we
identify the chosen $X_{\Delta}$'s with weighted projective planes or
quotients thereof by applying Lemma \ref{WEIGHTEDLEMMA}.

\section{Proof of Theorem \ref{THM3}: Step 1\label{STEP1}}

\begin{lemma}
\label{Lemma11}Among the $125$ possible combinations $(\alpha_{1},\alpha
_{2},\alpha_{3})$ of types of triples of pairs \emph{(\ref{TRIPLES}),} with
$\alpha_{1},\alpha_{2},\alpha_{3}\in\{\mathbf{1},\mathbf{2},\mathbf{3}%
,\mathbf{4},\mathbf{5}\},$ there are only $32$ satisfying simultaneously
conditions \emph{(\ref{MOD9I})} and \emph{(\ref{MOD9II});} namely,
\[
(\mathbf{1,3,4}),(\mathbf{1,4,5}),(\mathbf{2,3,4}),(\mathbf{2,4,5}%
),(\mathbf{3,3,3}),(\mathbf{3,3,5}),(\mathbf{3,5,5}),(\mathbf{5,5,5}),
\]
together with their permutations.\end{lemma}
{\tiny
\begin{table}[h!]%
\[
\begin{tabular}
[c]{|c|c|c|c|c|c|c|}\hline
\textbf{Case} & $\left[  \widehat{p}_{1}\right]  _{9}$ & $\left[
q_{1}\right]  _{9}$ & $\left[  p_{2}\right]  _{9}$ & $\left[  q_{2}\right]
_{9}$ & $\left[  \widehat{p}_{1}q_{2}+p_{2}q_{1}\right]  _{9}$ & $%
\begin{array}
[c]{c}%
(\ref{MOD9I})\text{ is true}\\
\text{only if}%
\end{array}
$\\\hline\hline
$\left(  \mathbf{1,1},\alpha_{3}\right)  $ & $2$ & $3$ & $2$ & $3$ & $3$ &
$\alpha_{3}=\mathbf{4}$\\\hline
$\left(  \mathbf{1,2},\alpha_{3}\right)  $ & $2$ & $3$ & $\in\left\{
2,5,8\right\}  $ & $3$ & $3$ & $\alpha_{3}=\mathbf{4}$\\\hline
$\left(  \mathbf{1,3},\alpha_{3}\right)  $ & $2$ & $3$ & $\in\left\{
1,4,7\right\}  $ & $0$ & $3$ & $\alpha_{3}=\mathbf{4}$\\\hline
$\left(  \mathbf{1,4},\alpha_{3}\right)  $ & $2$ & $3$ & $\in\left\{
2,5,8\right\}  $ & $6$ & $0$ & $\alpha_{3}\in\{\mathbf{3},\mathbf{5}%
\}$\\\hline
$\left(  \mathbf{1,5},\alpha_{3}\right)  $ & $2$ & $3$ & $\in\left\{
1,4,7\right\}  $ & $0$ & $3$ & $\alpha_{3}=\mathbf{4}$\\\hline
$\left(  \mathbf{2,1},\alpha_{3}\right)  $ & $\in\left\{  2,5,8\right\}  $ &
$3$ & $2$ & $3$ & $3$ & $\alpha_{3}=\mathbf{4}$\\\hline
$\left(  \mathbf{2,2},\alpha_{3}\right)  $ & $\in\left\{  2,5,8\right\}  $ &
$3$ & $\in\left\{  2,5,8\right\}  $ & $3$ & $3$ & $\alpha_{3}=\mathbf{4}%
$\\\hline
$\left(  \mathbf{2,3},\alpha_{3}\right)  $ & $\in\left\{  2,5,8\right\}  $ &
$3$ & $\in\left\{  1,4,7\right\}  $ & $0$ & $3$ & $\alpha_{3}=\mathbf{4}%
$\\\hline
$\left(  \mathbf{2,4},\alpha_{3}\right)  $ & $\in\left\{  2,5,8\right\}  $ &
$3$ & $\in\left\{  2,5,8\right\}  $ & $6$ & $0$ & $\alpha_{3}\in
\{\mathbf{3},\mathbf{5}\}$\\\hline
$\left(  \mathbf{2,5},\alpha_{3}\right)  $ & $\in\left\{  2,5,8\right\}  $ &
$3$ & $\in\left\{  1,4,7\right\}  $ & $0$ & $3$ & $\alpha_{3}=\mathbf{4}%
$\\\hline
$\left(  \mathbf{3,1},\alpha_{3}\right)  $ & $\in\left\{  1,4,7\right\}  $ &
$0$ & $2$ & $3$ & $3$ & $\alpha_{3}=\mathbf{4}$\\\hline
$\left(  \mathbf{3,2},\alpha_{3}\right)  $ & $\in\left\{  1,4,7\right\}  $ &
$0$ & $\in\left\{  2,5,8\right\}  $ & $3$ & $3$ & $\alpha_{3}=\mathbf{4}%
$\\\hline
$\left(  \mathbf{3,3},\alpha_{3}\right)  $ & $\in\left\{  1,4,7\right\}  $ &
$0$ & $\in\left\{  1,4,7\right\}  $ & $0$ & $0$ & $\alpha_{3}\in
\{\mathbf{3},\mathbf{5}\}$\\\hline
$\left(  \mathbf{3,4},\alpha_{3}\right)  $ & $\in\left\{  1,4,7\right\}  $ &
$0$ & $\in\left\{  2,5,8\right\}  $ & $6$ & $6$ & $\alpha_{3}\in
\{\mathbf{1},\mathbf{2}\}$\\\hline
$\left(  \mathbf{3,5},\alpha_{3}\right)  $ & $\in\left\{  1,4,7\right\}  $ &
$0$ & $\in\left\{  1,4,7\right\}  $ & $0$ & $0$ & $\alpha_{3}\in
\{\mathbf{3},\mathbf{5}\}$\\\hline
$\left(  \mathbf{4,1},\alpha_{3}\right)  $ & $\in\left\{  2,5,8\right\}  $ &
$6$ & $2$ & $3$ & $0$ & $\alpha_{3}\in\{\mathbf{3},\mathbf{5}\}$\\\hline
$\left(  \mathbf{4,2},\alpha_{3}\right)  $ & $\in\left\{  2,5,8\right\}  $ &
$6$ & $\in\left\{  2,5,8\right\}  $ & $3$ & $0$ & $\alpha_{3}\in
\{\mathbf{3},\mathbf{5}\}$\\\hline
$\left(  \mathbf{4,3},\alpha_{3}\right)  $ & $\in\left\{  2,5,8\right\}  $ &
$6$ & $\in\left\{  1,4,7\right\}  $ & $0$ & $6$ & $\alpha_{3}\in
\{\mathbf{1},\mathbf{2}\}$\\\hline
$\left(  \mathbf{4,4},\alpha_{3}\right)  $ & $\in\left\{  2,5,8\right\}  $ &
$6$ & $\in\left\{  2,5,8\right\}  $ & $6$ & $6$ & $\alpha_{3}\in
\{\mathbf{1},\mathbf{2}\}$\\\hline
$\left(  \mathbf{4,5},\alpha_{3}\right)  $ & $\in\left\{  2,5,8\right\}  $ &
$6$ & $\in\left\{  1,4,7\right\}  $ & $0$ & $6$ & $\alpha_{3}\in
\{\mathbf{1},\mathbf{2}\}$\\\hline
$\left(  \mathbf{5,1},\alpha_{3}\right)  $ & $\in\left\{  1,4,7\right\}  $ &
$0$ & $2$ & $3$ & $3$ & $\alpha_{3}=\mathbf{4}$\\\hline
$\left(  \mathbf{5,2},\alpha_{3}\right)  $ & $\in\left\{  1,4,7\right\}  $ &
$0$ & $\in\left\{  2,5,8\right\}  $ & $3$ & $3$ & $\alpha_{3}=\mathbf{4}%
$\\\hline
$\left(  \mathbf{5,3},\alpha_{3}\right)  $ & $\in\left\{  1,4,7\right\}  $ &
$0$ & $\in\left\{  1,4,7\right\}  $ & $0$ & $0$ & $\alpha_{3}\in
\{\mathbf{3},\mathbf{5}\}$\\\hline
$\left(  \mathbf{5,4},\alpha_{3}\right)  $ & $\in\left\{  1,4,7\right\}  $ &
$0$ & $\in\left\{  2,5,8\right\}  $ & $6$ & $6$ & $\alpha_{3}\in
\{\mathbf{1},\mathbf{2}\}$\\\hline
$\left(  \mathbf{5,5},\alpha_{3}\right)  $ & $\in\left\{  1,4,7\right\}  $ &
$0$ & $\in\left\{  1,4,7\right\}  $ & $0$ & $0$ & $\alpha_{3}\in
\{\mathbf{3},\mathbf{5}\}$\\\hline
\end{tabular}
\ \ \ %
\]
\smallskip\caption{{}}%
\label{Table2}%
\end{table}}
\vspace{-0.8cm}
\begin{table}[h!]
{\tiny
\[%
\begin{tabular}
[c]{|c|c|c|c|c|c|c|}\hline
\textbf{Case} & $\left[  p_{1}\right]  _{9}$ & $\left[  q_{1}\right]  _{9}$ &
$\left[  \widehat{p}_{3}\right]  _{9}$ & $\left[  q_{3}\right]  _{9}$ &
$\left[  p_{1}q_{3}+\widehat{p}_{3}q_{1}\right]  _{9}$ & $%
\begin{array}
[c]{c}%
(\ref{MOD9II})\text{ is true}\\
\text{only if}%
\end{array}
$\\\hline\hline
$\left(  \mathbf{1,}\alpha_{2}\mathbf{,1}\right)  $ & $2$ & $3$ & $2$ & $3$ &
$3$ & $\alpha_{2}=\mathbf{4}$\\\hline
$\left(  \mathbf{1,}\alpha_{2}\mathbf{,2}\right)  $ & $2$ & $3$ & $\in\left\{
2,5,8\right\}  $ & $3$ & $3$ & $\alpha_{2}=\mathbf{4}$\\\hline
$\left(  \mathbf{1,}\alpha_{2}\mathbf{,3}\right)  $ & $2$ & $3$ & $\in\left\{
1,4,7\right\}  $ & $0$ & $3$ & $\alpha_{2}=\mathbf{4}$\\\hline
$\left(  \mathbf{1,}\alpha_{2}\mathbf{,4}\right)  $ & $2$ & $3$ & $\in\left\{
2,5,8\right\}  $ & $6$ & $0$ & $\alpha_{2}\in\{\mathbf{3},\mathbf{5}%
\}$\\\hline
$\left(  \mathbf{1,}\alpha_{2}\mathbf{,5}\right)  $ & $2$ & $3$ & $\in\left\{
1,4,7\right\}  $ & $0$ & $3$ & $\alpha_{2}=\mathbf{4}$\\\hline
$\left(  \mathbf{2,}\alpha_{2}\mathbf{,1}\right)  $ & $\in\left\{
2,5,8\right\}  $ & $3$ & $2$ & $3$ & $3$ & $\alpha_{2}=\mathbf{4}$\\\hline
$\left(  \mathbf{2,}\alpha_{2}\mathbf{,2}\right)  $ & $\in\left\{
2,5,8\right\}  $ & $3$ & $\in\left\{  2,5,8\right\}  $ & $3$ & $3$ &
$\alpha_{2}=\mathbf{4}$\\\hline
$\left(  \mathbf{2,}\alpha_{2}\mathbf{,3}\right)  $ & $\in\left\{
2,5,8\right\}  $ & $3$ & $\in\left\{  1,4,7\right\}  $ & $0$ & $3$ &
$\alpha_{2}=\mathbf{4}$\\\hline
$\left(  \mathbf{2,}\alpha_{2}\mathbf{,4}\right)  $ & $\in\left\{
2,5,8\right\}  $ & $3$ & $\in\left\{  2,5,8\right\}  $ & $6$ & $0$ &
$\alpha_{2}\in\{\mathbf{3},\mathbf{5}\}$\\\hline
$\left(  \mathbf{2,}\alpha_{2}\mathbf{,5}\right)  $ & $\in\left\{
2,5,8\right\}  $ & $3$ & $\in\left\{  1,4,7\right\}  $ & $0$ & $3$ &
$\alpha_{2}=\mathbf{4}$\\\hline
$\left(  \mathbf{3,}\alpha_{2}\mathbf{,1}\right)  $ & $\in\left\{
1,4,7\right\}  $ & $0$ & $2$ & $3$ & $3$ & $\alpha_{2}=\mathbf{4}$\\\hline
$\left(  \mathbf{3,}\alpha_{2}\mathbf{,2}\right)  $ & $\in\left\{
1,4,7\right\}  $ & $0$ & $\in\left\{  2,5,8\right\}  $ & $3$ & $3$ &
$\alpha_{2}=\mathbf{4}$\\\hline
$\left(  \mathbf{3,}\alpha_{2}\mathbf{,3}\right)  $ & $\in\left\{
1,4,7\right\}  $ & $0$ & $\in\left\{  1,4,7\right\}  $ & $0$ & $0$ &
$\alpha_{2}\in\{\mathbf{3},\mathbf{5}\}$\\\hline
$\left(  \mathbf{3,}\alpha_{2}\mathbf{,4}\right)  $ & $\in\left\{
1,4,7\right\}  $ & $0$ & $\in\left\{  2,5,8\right\}  $ & $6$ & $6$ &
$\alpha_{2}\in\{\mathbf{1},\mathbf{2}\}$\\\hline
$\left(  \mathbf{3,}\alpha_{2}\mathbf{,5}\right)  $ & $\in\left\{
1,4,7\right\}  $ & $0$ & $\in\left\{  1,4,7\right\}  $ & $0$ & $0$ &
$\alpha_{2}\in\{\mathbf{3},\mathbf{5}\}$\\\hline
$\left(  \mathbf{4,}\alpha_{2}\mathbf{,1}\right)  $ & $\in\left\{
2,5,8\right\}  $ & $6$ & $2$ & $3$ & $0$ & $\alpha_{2}\in\{\mathbf{3}%
,\mathbf{5}\}$\\\hline
$\left(  \mathbf{4,}\alpha_{2}\mathbf{,2}\right)  $ & $\in\left\{
2,5,8\right\}  $ & $6$ & $\in\left\{  2,5,8\right\}  $ & $3$ & $0$ &
$\alpha_{2}\in\{\mathbf{3},\mathbf{5}\}$\\\hline
$\left(  \mathbf{4,}\alpha_{2}\mathbf{,3}\right)  $ & $\in\left\{
2,5,8\right\}  $ & $6$ & $\in\left\{  1,4,7\right\}  $ & $0$ & $6$ &
$\alpha_{2}\in\{\mathbf{1},\mathbf{2}\}$\\\hline
$\left(  \mathbf{4,}\alpha_{2}\mathbf{,4}\right)  $ & $\in\left\{
2,5,8\right\}  $ & $6$ & $\in\left\{  2,5,8\right\}  $ & $6$ & $6$ &
$\alpha_{2}\in\{\mathbf{1},\mathbf{2}\}$\\\hline
$\left(  \mathbf{4,}\alpha_{2}\mathbf{,5}\right)  $ & $\in\left\{
2,5,8\right\}  $ & $6$ & $\in\left\{  1,4,7\right\}  $ & $0$ & $6$ &
$\alpha_{2}\in\{\mathbf{1},\mathbf{2}\}$\\\hline
$\left(  \mathbf{5,}\alpha_{2}\mathbf{,1}\right)  $ & $\in\left\{
1,4,7\right\}  $ & $0$ & $2$ & $3$ & $3$ & $\alpha_{2}=\mathbf{4}$\\\hline
$\left(  \mathbf{5,}\alpha_{2}\mathbf{,2}\right)  $ & $\in\left\{
1,4,7\right\}  $ & $0$ & $\in\left\{  2,5,8\right\}  $ & $3$ & $3$ &
$\alpha_{2}=\mathbf{4}$\\\hline
$\left(  \mathbf{5,}\alpha_{2}\mathbf{,3}\right)  $ & $\in\left\{
1,4,7\right\}  $ & $0$ & $\in\left\{  1,4,7\right\}  $ & $0$ & $0$ &
$\alpha_{2}\in\{\mathbf{3},\mathbf{5}\}$\\\hline
$\left(  \mathbf{5,}\alpha_{2}\mathbf{,4}\right)  $ & $\in\left\{
1,4,7\right\}  $ & $0$ & $\in\left\{  2,5,8\right\}  $ & $6$ & $6$ &
$\alpha_{2}\in\{\mathbf{1},\mathbf{2}\}$\\\hline
$\left(  \mathbf{5,}\alpha_{2}\mathbf{,5}\right)  $ & $\in\left\{
1,4,7\right\}  $ & $0$ & $\in\left\{  1,4,7\right\}  $ & $0$ & $0$ &
$\alpha_{2}\in\{\mathbf{3},\mathbf{5}\}$\\\hline
\end{tabular}
\ \ \ \ \
\]
}\smallskip\caption{{}}%
\label{Table3}%
\end{table}
\newpage

\begin{proof}
By Table \ref{Table2} there are $38$ combinations $(\alpha_{1},\alpha
_{2},\alpha_{3})$ of types of triples (\ref{TRIPLES}), with $\alpha_{1}%
,\alpha_{2},\alpha_{3}\in\{\mathbf{1},\mathbf{2},\mathbf{3},\mathbf{4}%
,\mathbf{5}\},$ satisfying condition (\ref{MOD9I}). Correspondingly, Table \ref{Table3} shows that there are
$38$ combinations $(\alpha_{1},\alpha_{2},\alpha_{3})$ of types of triples
(\ref{TRIPLES}), with $\alpha_{1},\alpha_{2},\alpha_{3}\in\{\mathbf{1}%
,\mathbf{2},\mathbf{3},\mathbf{4},\mathbf{5}\},$ satisfying condition
(\ref{MOD9II}). Obviously, the combinations $(\alpha_{1}%
,\alpha_{2},\alpha_{3})$ of types of triples of pairs (\ref{TRIPLES}), with
$\alpha_{1},\alpha_{2},\alpha_{3}\in\{\mathbf{1},\mathbf{2},\mathbf{3}%
,\mathbf{4},\mathbf{5}\},$ satisfying \textit{both} (\ref{MOD9I}) and
(\ref{MOD9II}), are the $32$ combinations given in the statement of Lemma.
\end{proof}

\begin{lemma}
\label{Lemma12} There are no admissible triples of pairs \emph{(\ref{TRIPLES})
}among those corresponding to the $125$ type combinations $(\alpha_{1}%
,\alpha_{2},\alpha_{3})$ with $\alpha_{1},\alpha_{2},\alpha_{3}\in
\{\mathbf{1},\mathbf{2},\mathbf{3},\mathbf{4},\mathbf{5}\}.$
\end{lemma}

\noindent{}\textit{Sketch of proof}. First, we express the triples of pairs
$\left\{  \left.  (p_{i},q_{i})\in\mathbb{Z}^{2}\right\vert 1\leq
i\leq3\right\}  $ corresponding to the $32$ type combinations $(\alpha
_{1},\alpha_{2},\alpha_{3})$ found in Lemma \ref{Lemma11} in terms of $\xi
_{i}$ for $i\in\{1,2,3\}$ as in Table \ref{Table1}. Setting%
\[
\mathfrak{A}_{j}:=\left\{  \left(  \xi_{1},\xi_{2},\xi_{3}\right)
\in\mathbb{Z}^{3}\left\vert
\begin{array}
[c]{c}%
\xi_{1}+\xi_{2}+\xi_{3}\leq10,\ \xi_{j}\geq0,\smallskip\text{ }\\
\text{and }\xi_{k}\geq1,\ \forall k\in\{1,2,3\}\mathbb{r}\{j\}
\end{array}
\right.  \right\}  ,
\]
for $j\in\{1,2,3\},$
\[
\mathfrak{A}_{j,k}:=\left\{  \left(  \xi_{1},\xi_{2},\xi_{3}\right)
\in\mathbb{Z}^{3}\left\vert
\begin{array}
[c]{c}%
\xi_{1}+\xi_{2}+\xi_{3}\leq10,\ \xi_{j}=0,\xi_{k}\geq0,\smallskip\text{ }\\
\text{and }\xi_{\mu}\geq1,\ \forall\mu\in\{1,2,3\}\mathbb{r}\{j,k\}
\end{array}
\right.  \right\}  ,
\]
for $j,k\in\{1,2,3\},$ $j\neq k,$ and
\[
\mathfrak{B}:=\left\{  \left.  \left(  \xi_{1},\xi_{2},\xi_{3}\right)
\in\mathbb{Z}^{3}\right\vert \xi_{1}+\xi_{2}+\xi_{3}\leq11,\ \xi_{1},\xi
_{2},\xi_{3}\geq1\right\}  ,
\]
we explain what condition (\ref{Condition3}) means for each of these $32$
cases in Table \ref{Table4}.
\begin{table}[h]
{\small
\[%
\begin{tabular}
[c]{|c|c|c|c|}\hline
\textbf{Case} & Condition (\ref{Condition3}) & \textbf{Case} & Condition
(\ref{Condition3})\\\hline\hline
$(\mathbf{1,3,4})$ & $\left(  \xi_{1},\xi_{2},\xi_{3}\right)  \in
\mathfrak{A}_{1,3}$ & $(\mathbf{4,1,3})$ & $\left(  \xi_{1},\xi_{2},\xi
_{3}\right)  \in\mathfrak{A}_{2,1}$\\\hline
$(\mathbf{1,4,3})$ & $\left(  \xi_{1},\xi_{2},\xi_{3}\right)  \in
\mathfrak{A}_{1,2}$ & $(\mathbf{4,1,5})$ & $\left(  \xi_{1},\xi_{2},\xi
_{3}\right)  \in\mathfrak{A}_{2,1}$\\\hline
$(\mathbf{1,4,5})$ & $\left(  \xi_{1},\xi_{2},\xi_{3}\right)  \in
\mathfrak{A}_{1,2}$ & $(\mathbf{4,2,3})$ & $\left(  \xi_{1},\xi_{2},\xi
_{3}\right)  \in\mathfrak{A}_{1}$\\\hline
$(\mathbf{1,5,4})$ & $\left(  \xi_{1},\xi_{2},\xi_{3}\right)  \in
\mathfrak{A}_{1,3}$ & $(\mathbf{4,2,5})$ & $\left(  \xi_{1},\xi_{2},\xi
_{3}\right)  \in\mathfrak{A}_{1}$\\\hline
$(\mathbf{2,3,4})$ & $\left(  \xi_{1},\xi_{2},\xi_{3}\right)  \in
\mathfrak{A}_{3}$ & $(\mathbf{4,3,1})$ & $\left(  \xi_{1},\xi_{2},\xi
_{3}\right)  \in\mathfrak{A}_{3,1}$\\\hline
$(\mathbf{2,4,3})$ & $\left(  \xi_{1},\xi_{2},\xi_{3}\right)  \in
\mathfrak{A}_{2}$ & $(\mathbf{4,3,2})$ & $\left(  \xi_{1},\xi_{2},\xi
_{3}\right)  \in\mathfrak{A}_{1}$\\\hline
$(\mathbf{2,4,5})$ & $\left(  \xi_{1},\xi_{2},\xi_{3}\right)  \in
\mathfrak{A}_{2}$ & $(\mathbf{4,5,1})$ & $\left(  \xi_{1},\xi_{2},\xi
_{3}\right)  \in\mathfrak{A}_{3,1}$\\\hline
$(\mathbf{2,5,4})$ & $\left(  \xi_{1},\xi_{2},\xi_{3}\right)  \in
\mathfrak{A}_{3}$ & $(\mathbf{4,5,2})$ & $\left(  \xi_{1},\xi_{2},\xi
_{3}\right)  \in\mathfrak{A}_{1}$\\\hline
$(\mathbf{3,1,4})$ & $\left(  \xi_{1},\xi_{2},\xi_{3}\right)  \in
\mathfrak{A}_{2,3}$ & $(\mathbf{5,1,4})$ & $\left(  \xi_{1},\xi_{2},\xi
_{3}\right)  \in\mathfrak{A}_{2,3}$\\\hline
$(\mathbf{3,2,4})$ & $\left(  \xi_{1},\xi_{2},\xi_{3}\right)  \in
\mathfrak{A}_{3}$ & $(\mathbf{5,2,4})$ & $\left(  \xi_{1},\xi_{2},\xi
_{3}\right)  \in\mathfrak{A}_{3}$\\\hline
$(\mathbf{3,3,3})$ & $\left(  \xi_{1},\xi_{2},\xi_{3}\right)  \in\mathfrak{B}$
& $(\mathbf{5,3,3})$ & $\left(  \xi_{1},\xi_{2},\xi_{3}\right)  \in
\mathfrak{B}$\\\hline
$(\mathbf{3,3,5})$ & $\left(  \xi_{1},\xi_{2},\xi_{3}\right)  \in\mathfrak{B}$
& $(\mathbf{5,3,5})$ & $\left(  \xi_{1},\xi_{2},\xi_{3}\right)  \in
\mathfrak{B}$\\\hline
$(\mathbf{3,4,1})$ & $\left(  \xi_{1},\xi_{2},\xi_{3}\right)  \in
\mathfrak{A}_{3,2}$ & $(\mathbf{5,4,1})$ & $\left(  \xi_{1},\xi_{2},\xi
_{3}\right)  \in\mathfrak{A}_{3,2}$\\\hline
$(\mathbf{3,4,2})$ & $\left(  \xi_{1},\xi_{2},\xi_{3}\right)  \in
\mathfrak{A}_{2}$ & $(\mathbf{5,4,2})$ & $\left(  \xi_{1},\xi_{2},\xi
_{3}\right)  \in\mathfrak{A}_{2}$\\\hline
$(\mathbf{3,5,3})$ & $\left(  \xi_{1},\xi_{2},\xi_{3}\right)  \in\mathfrak{B}$
& $(\mathbf{5,5,3})$ & $\left(  \xi_{1},\xi_{2},\xi_{3}\right)  \in
\mathfrak{B}$\\\hline
$(\mathbf{3,5,5})$ & $\left(  \xi_{1},\xi_{2},\xi_{3}\right)  \in\mathfrak{B}$
& $(\mathbf{5,5,5})$ & $\left(  \xi_{1},\xi_{2},\xi_{3}\right)  \in
\mathfrak{B}$\\\hline
\end{tabular}
\ \ \ \ \ \ \
\]
} \smallskip\caption{{}}%
\label{Table4}%
\end{table}
\vspace{-0.2cm}

\noindent Note that
\[
\sharp(\mathfrak{A}_{j})=\sum_{\kappa=2}^{10}\tbinom{\kappa-1}{1}+\sum
_{\kappa=3}^{10}\tbinom{\kappa-1}{2}=\allowbreak165,\ \sharp(\mathfrak{A}%
_{j,k})=55,\ \sharp(\mathfrak{B})=\tbinom{11}{3}=165.
\]
One can, of course, test directly the validity of (\ref{Condition1}) and
(\ref{Condition2}) for all these possibilities. Nevertheless, there is a more
economic way to proceed by using reductio ad absurdum. Let us discuss it
exemplarily in the case $(\mathbf{2,3,4})$ in which{\small
\[%
\begin{tabular}
[c]{|c|c|c|c|c|c|c|c|c|}\hline
$p_{1}$ & $\widehat{p}_{1}$ & $q_{1}$ & $p_{2}$ & $\widehat{p}_{2}$ & $q_{2}$
& $p_{3}$ & $\widehat{p}_{3}$ & $q_{3}$\\\hline\hline
$3\xi_{1}+2$ & $3\xi_{1}+2$ & $9\xi_{1}+3$ & $3\xi_{2}+1$ & $6\xi_{2}+1$ &
$9\xi_{2}$ & $6\xi_{3}+5$ & $6\xi_{3}+5$ & $9\xi_{3}+6$\\\hline
\end{tabular}
\ \ \ \ \
\]
}for a $3$-tuple $\left(  \xi_{1},\xi_{2},\xi_{3}\right)  \in\mathfrak{A}%
_{3}.$ If $\left\{  \left.  (p_{i},q_{i})\in\mathbb{Z}^{2}\right\vert 1\leq
i\leq3\right\}  $ were an admissible triple of pairs, then (\ref{Condition1})
would give%
\[
\left(  9\xi_{1}+3\right)  \left(  9\xi_{2}\right)  =\allowbreak27\xi
_{2}+81\xi_{1}\xi_{2}\mid\allowbreak9\xi_{1}+27\xi_{2}+9\xi_{3}+54\xi_{1}%
\xi_{2}+9,\text{ i.e.,}%
\]%
\begin{equation}
3\left(  3\xi_{1}+1\right)  \xi_{2}\mid\xi_{1}+3\xi_{2}+\xi_{3}+6\xi_{1}%
\xi_{2}+1=3\left(  3\xi_{1}+1\right)  \xi_{2}-3\xi_{1}\xi_{2}+\xi_{1}+\xi
_{3}+1, \label{test1}%
\end{equation}
meaning that $3\left(  3\xi_{1}+1\right)  \xi_{2}\mid3\xi_{1}\xi_{2}-\xi
_{1}-\xi_{3}-1.$ Therefore, $3\xi_{1}\xi_{2}-\xi_{1}-\xi_{3}-1\leq0$ (because
otherwise we would deduce that $3\xi_{2}+9\xi_{1}\xi_{2}\leq3\xi_{1}\xi
_{2}-\xi_{1}-\xi_{3}-1,$ i.e., that $10\leq\xi_{1}+3\xi_{2}+\xi_{3}+6\xi
_{1}\xi_{2}\leq-1,$ a contradiction). Consequently,%
\begin{align}
3\xi_{1}\xi_{2}-1  &  \leq\xi_{1}+\xi_{3}\leq10-\xi_{2}\Longrightarrow
4\leq(3\xi_{1}+1)\xi_{2}\leq11\nonumber\\
&  \Longrightarrow\left(  \xi_{1},\xi_{2}\right)  \in\left\{
(1,1),(1,2),(2,1),(3,1)\right\}  . \label{test_2}%
\end{align}
Since $0\leq\xi_{3}\leq8,$ (\ref{test1}) and (\ref{test_2})\ would determine
the values of $\xi_{3}$ as follows: {\small
\[%
\begin{tabular}
[c]{|c|c|c|c|c|c|c|c|}\hline
$\left(  \xi_{1},\xi_{2},\xi_{3}\right)  $ & $p_{1}$ & $q_{1}$ & $q_{2}$ &
$\widehat{p}_{3}$ & $q_{3}$ & $q_{1}q_{3}$ & $p_{1}q_{3}+\widehat{p}_{3}%
q_{1}+q_{2}$\\\hline\hline
$(1,1,1)$ & $5$ & $12$ & $9$ & $11$ & $15$ & $180$ & $216$\\\hline
$(1,2,4)$ & $5$ & $12$ & $18$ & $29$ & $42$ & $504$ & $576$\\\hline
$(2,1,3)$ & $8$ & $21$ & $9$ & $23$ & $33$ & $693$ & $756$\\\hline
$(3,1,5)$ & $11$ & $30$ & $9$ & $35$ & $51$ & $1530$ & $1620$\\\hline
\end{tabular}
\ \ \ \ \
\]
}Hence, these four $3$-tuples $\left(  \xi_{1},\xi_{2},\xi_{3}\right)
\in\mathfrak{A}_{3}$ would provide numbers $p_{1},q_{1},$ $q_{2},$
$\widehat{p}_{3},$ $q_{3}$ which do not satisfy (\ref{Condition2})! Using
analogous arguments one shows that none of the remaining $31$ cases leads to
admissible triples of pairs.\hfill$\square$

\section{Proof of Theorem \ref{THM3}: Step 2\label{STEP2}}

\begin{lemma}
\label{Lemma21}There are no admissible triples of pairs \emph{(\ref{TRIPLES})
}among those corresponding to the type combinations $(\alpha_{1},\alpha
_{2},\alpha_{3})$ with $\alpha_{1}\in\{\mathbf{1},\mathbf{2},\mathbf{3}%
,\mathbf{4},\mathbf{5}\}$ and
\[
(\alpha_{2},\alpha_{3})\in\left(  \{\mathbf{1},\mathbf{2},\mathbf{3}%
,\mathbf{4},\mathbf{5}\}\times\{\mathbf{7}\}\right)  \cup\left(
\{\mathbf{7}\}\times\{\mathbf{1},\mathbf{2},\mathbf{3},\mathbf{4}%
,\mathbf{5}\}\right)  .
\]

\end{lemma}

\begin{proof}
If $\alpha_{1},\alpha_{2}\in\{\mathbf{1},\mathbf{2},\mathbf{3},\mathbf{4}%
,\mathbf{5}\}$ and $\alpha_{3}=\mathbf{7},$ then $\left[  \widehat{p}_{1}%
q_{2}+p_{2}q_{1}\right]  _{9}\in\{0,3,6\}$ (cf. the sixth column of Table
\ref{Table2}) and $q_{3}=1,$ i.e., $\left[  \widehat{p}_{1}q_{2}+p_{2}%
q_{1}+q_{3}\right]  _{9}\in\{1,4,7\}.$ Thus, condition (\ref{MOD9I}) is not
satisfied. Analogously, one shows that condition (\ref{MOD9II}) is not
satisfied whenever $\alpha_{1},\alpha_{3}\in\{\mathbf{1},\mathbf{2}%
,\mathbf{3},\mathbf{4},\mathbf{5}\}$ and $\alpha_{2}=\mathbf{7}.$
\end{proof}

\begin{lemma}
\label{Lemma22}There exist exactly $10$ admissible triples of pairs
\emph{(\ref{TRIPLES}) }among those corresponding to the type combinations
$(\alpha_{1},\alpha_{2},\alpha_{3})$ with $\alpha_{1}\in\{\mathbf{1}%
,\mathbf{2},\mathbf{3},\mathbf{4},\mathbf{5}\}$ and
\[
(\alpha_{2},\alpha_{3})\in\left(  \{\mathbf{1},\mathbf{2},\mathbf{3}%
,\mathbf{4},\mathbf{5}\}\times\{\mathbf{6}\}\right)  \cup\left(
\{\mathbf{6}\}\times\{\mathbf{1},\mathbf{2},\mathbf{3},\mathbf{4}%
,\mathbf{5}\}\right)  .
\]

\end{lemma}

\noindent{}\textit{Sketch of proof}. For $\alpha_{1},\alpha_{2}\in
\{\mathbf{1},\mathbf{2},\mathbf{3},\mathbf{4},\mathbf{5}\}$ and $\alpha
_{3}=\mathbf{6}$ we build Table \ref{Table5}. In its second column we tabulate
$\left[  \widehat{p}_{1}q_{2}+p_{2}q_{1}\right]  _{9}$ (cf. the sixth column
of Table \ref{Table2}). After having expressed $q_{1},$ $q_{2}$ in terms of
$\xi_{1},\xi_{2}$ (as in Table \ref{Table1})\ we write the restrictions
(inequalities) coming from (\ref{Condition3}) in its third column. The fourth
column contains the values of $q_{3}$ so that both (\ref{Condition3}) and
(\ref{MOD9I}) are true. (In particular, in the case $\left(  \mathbf{2,2}%
,\mathbf{6}\right)  $ the expected value $q_{3}=6$ is impossible because
$\xi_{1},\xi_{2}\geq1.$) Finally, the last column informs us whether
(\ref{MOD9II}) \ is true for these $q_{3}$'s.
\newpage

\begin{table}[h]
{\tiny
\[%
\begin{tabular}
[c]{|c|c|c|c|c|}\hline
\textbf{Case} & $\left[  \widehat{p}_{1}q_{2}+p_{2}q_{1}\right]  _{9}$ & $%
\begin{array}
[c]{c}%
\text{(\ref{Condition3}) is true}\\
\text{whenever}%
\end{array}
$ & $%
\begin{array}
[c]{c}%
\text{(\ref{Condition3}) \& (\ref{MOD9I}) true}\\
\text{only if }q_{3}\text{ equals}%
\end{array}
$ & $%
\begin{array}
[c]{c}%
\text{Is (\ref{MOD9II}) true}\\
\text{for these }q_{3}\text{'s?}%
\end{array}
$\\\hline\hline
$\left(  \mathbf{1,1},\mathbf{6}\right)  $ & $3$ & $2\leq q_{3}\leq7$ & $6$ &
YES\\\hline
$\left(  \mathbf{1,2},\mathbf{6}\right)  $ & $3$ & $3\leq\xi_{2}+q_{3}\leq7$ &
$6$ & YES\\\hline
$\left(  \mathbf{1,3},\mathbf{6}\right)  $ & $3$ & $3\leq\xi_{2}+q_{3}\leq9$ &
$6$ & NO\\\hline
$\left(  \mathbf{1,4},\mathbf{6}\right)  $ & $0$ & $2\leq\xi_{2}+q_{3}\leq10$
& $9$ & YES\\\hline
$\left(  \mathbf{1,5},\mathbf{6}\right)  $ & $3$ & $3\leq\xi_{2}+q_{3}\leq9$ &
$6$ & NO\\\hline
$\left(  \mathbf{2,1},\mathbf{6}\right)  $ & $3$ & $3\leq\xi_{1}+q_{3}\leq7$ &
$6$ & YES\\\hline
$\left(  \mathbf{2,2},\mathbf{6}\right)  $ & $3$ & $4\leq\xi_{1}+\xi_{2}%
+q_{3}\leq7$ & $6$ (impossible) & -----\\\hline
$\left(  \mathbf{2,3},\mathbf{6}\right)  $ & $3$ & $4\leq\xi_{1}+\xi_{2}%
+q_{3}\leq9$ & $6$ & NO\\\hline
$\left(  \mathbf{2,4},\mathbf{6}\right)  $ & $0$ & $3\leq\xi_{1}+\xi_{2}%
+q_{3}\leq10$ & $9$ & YES\\\hline
$\left(  \mathbf{2,5},\mathbf{6}\right)  $ & $3$ & $4\leq\xi_{1}+\xi_{2}%
+q_{3}\leq9$ & $6$ & NO\\\hline
$\left(  \mathbf{3,1},\mathbf{6}\right)  $ & $3$ & $3\leq\xi_{1}+q_{3}\leq9$ &
$6$ & NO\\\hline
$\left(  \mathbf{3,2},\mathbf{6}\right)  $ & $3$ & $4\leq\xi_{1}+\xi_{2}%
+q_{3}\leq9$ & $6$ & NO\\\hline
$\left(  \mathbf{3,3},\mathbf{6}\right)  $ & $0$ & $4\leq\xi_{1}+\xi_{2}%
+q_{3}\leq11$ & $9$ & YES\\\hline
$\left(  \mathbf{3,4},\mathbf{6}\right)  $ & $6$ & $3\leq\xi_{1}+\xi_{2}%
+q_{3}\leq11$ & $3$ & YES\\\hline
$\left(  \mathbf{3,5},\mathbf{6}\right)  $ & $0$ & $4\leq\xi_{1}+\xi_{2}%
+q_{3}\leq11$ & $9$ & YES\\\hline
$\left(  \mathbf{4,1},\mathbf{6}\right)  $ & $0$ & $2\leq\xi_{1}+q_{3}\leq10$
& $9$ & YES\\\hline
$\left(  \mathbf{4,2},\mathbf{6}\right)  $ & $0$ & $3\leq\xi_{1}+\xi_{2}%
+q_{3}\leq10$ & $9$ & YES\\\hline
$\left(  \mathbf{4,3},\mathbf{6}\right)  $ & $6$ & $3\leq\xi_{1}+\xi_{2}%
+q_{3}\leq11$ & $3$ & NO\\\hline
$\left(  \mathbf{4,4},\mathbf{6}\right)  $ & $6$ & $2\leq\xi_{1}+\xi_{2}%
+q_{3}\leq12$ & $3$ or $12$ & YES\\\hline
$\left(  \mathbf{4,5},\mathbf{6}\right)  $ & $6$ & $3\leq\xi_{1}+\xi_{2}%
+q_{3}\leq11$ & $3$ & NO\\\hline
$\left(  \mathbf{5,1},\mathbf{6}\right)  $ & $3$ & $3\leq\xi_{1}+q_{3}\leq9$ &
$6$ & NO\\\hline
$\left(  \mathbf{5,2},\mathbf{6}\right)  $ & $3$ & $4\leq\xi_{1}+\xi_{2}%
+q_{3}\leq9$ & $6$ & NO\\\hline
$\left(  \mathbf{5,3},\mathbf{6}\right)  $ & $0$ & $4\leq\xi_{1}+\xi_{2}%
+q_{3}\leq11$ & $9$ & YES\\\hline
$\left(  \mathbf{5,4},\mathbf{6}\right)  $ & $6$ & $3\leq\xi_{1}+\xi_{2}%
+q_{3}\leq11$ & $3$ & YES\\\hline
$\left(  \mathbf{5,5},\mathbf{6}\right)  $ & $0$ & $4\leq\xi_{1}+\xi_{2}%
+q_{3}\leq11$ & $9$ & YES\\\hline
\end{tabular}
\ \
\]
}\smallskip\caption{{}}%
\label{Table5}%
\end{table}
\vspace{-0.3cm}

\noindent Next, we analyze in detail the $14$ cases for which the
answer is \textsc{\textquotedblleft yes\textquotedblright}.\smallskip

\noindent{}$\bullet$ In the case $(\mathbf{1,1,6})$ we have $q_{3}=6$ and we
obtain just one admissible triple of pairs:%
\begin{equation}%
\begin{tabular}
[c]{|c|c|c|c|c|c|}\hline
$p_{1}$ & $q_{1}$ & $p_{2}$ & $q_{2}$ & $p_{3}$ & $q_{3}$\\\hline\hline
$2$ & $3$ & $2$ & $3$ & $1$ & $6$\\\hline
\end{tabular}
\ \ \label{AB1}%
\end{equation}

\noindent{}$\bullet$ In cases $(\mathbf{1,2,6})$ and $(\mathbf{2,1,6})$ we
have $\xi_{2}=1,q_{3}=6,$ and $\xi_{1}=1,q_{3}=6,$ respectively, and
(\ref{Condition1}) cannot be satisfied (because $36\nmid45$). Hence, there are
no admissible triples of pairs.\smallskip\

\noindent{}$\bullet$ In cases $(\mathbf{1,4,6})$ and $(\mathbf{4,1,6})$ we
have $\xi_{2}\in\{0,1\},q_{3}=9,$ and $\xi_{1}\in\{0,1\},q_{3}=9,$
respectively, and (\ref{Condition1}) cannot be satisfied for $\xi_{2}=1,$
resp. for $\xi_{1}=1$ (because $45\nmid72$). For this reason, the only triples
of pairs which are admissible (i.e., for which both (\ref{Condition1}) and
(\ref{Condition2}) are satified) are%
\begin{equation}%
\begin{tabular}
[c]{|c|c|c|c|c|c|}\hline
$p_{1}$ & $q_{1}$ & $p_{2}$ & $q_{2}$ & $p_{3}$ & $q_{3}$\\\hline\hline
$2$ & $3$ & $5$ & $6$ & $1$ & $9$\\\hline
\end{tabular}
\ \ \label{Ex1}%
\end{equation}
and%
\begin{equation}%
\begin{tabular}
[c]{|c|c|c|c|c|c|}\hline
$p_{1}$ & $q_{1}$ & $p_{2}$ & $q_{2}$ & $p_{3}$ & $q_{3}$\\\hline\hline
$5$ & $6$ & $2$ & $3$ & $1$ & $9$\\\hline
\end{tabular}
\ \ \label{Ex2}%
\end{equation}
$\bullet$ In cases $(\mathbf{2,4,6})$ and $(\mathbf{4,2,6})$ we have
necessarily $\xi_{1}=1,$ $\xi_{2}=0,$ $q_{3}=9,$ and $\xi_{1}=0,$ $\xi_{2}=1,$
$q_{3}=9,$ respectively, and (\ref{Condition1}) cannot be satisfied (because
$72\nmid99$). Hence, there are no admissible triples of pairs.\smallskip\

\noindent{}$\bullet$ In cases $\left(  \mathbf{3,3},\mathbf{6}\right)  $ and
$\left(  \mathbf{5,5},\mathbf{6}\right)  $ we have necessarily $\xi_{1}%
=\xi_{2}=1,$ $q_{3}=9,$ and (\ref{Condition1}) cannot be satisfied (because
$81\nmid108$). Therefore, there are no admissible triples of pairs.\smallskip\

\noindent{}$\bullet$ In cases $\left(  \mathbf{3,4},\mathbf{6}\right)  $ and
$\left(  \mathbf{5,4},\mathbf{6}\right)  $ we have $q_{3}=3$ and $\xi_{1}%
+\xi_{2}\in\{1,\ldots,8\}$ with $\xi_{1}\geq1$ and $\xi_{2}\geq0.$ If
(\ref{Condition1}) were true, then in particular $q_{1}\mid\widehat{p}%
_{1}q_{2}+q_{3},$ i.e.,%
\[
\xi_{1}\mid\xi_{2}+1\Longrightarrow(\xi_{1},\xi_{2})\in\{\left.
(1,j)\right\vert 0\leq j\leq7\}\cup\left\{
(2,1),(2,3),(2,5),(3,2),(3,5),(4,3)\right\}  .
\]
As for everyone of the $14$ possible values of $(\xi_{1},\xi_{2})$ at least
one of the divisibility conditions (\ref{Condition1}) and (\ref{Condition2})
is violated$,$ there are no admissible triples of pairs.\smallskip\

\noindent{}$\bullet$ Case $\left(  \mathbf{3,5},\mathbf{6}\right)  $: $\xi
_{1}=\xi_{2}=1,$ $q_{3}=9,$ and (\ref{Condition1}) cannot be satisfied
(because $81\nmid135$); no admissible triples of pairs occur.\smallskip\

\noindent{}$\bullet$ Case $\left(  \mathbf{4,4},\mathbf{6}\right)  $: Here,
\textit{either} $\xi_{1}=\xi_{2}=0,$ $q_{3}=12,$ giving the admissible triple
of pairs:%
\begin{equation}%
\begin{tabular}
[c]{|c|c|c|c|c|c|}\hline
$p_{1}$ & $q_{1}$ & $p_{2}$ & $q_{2}$ & $p_{3}$ & $q_{3}$\\\hline\hline
$5$ & $6$ & $5$ & $6$ & $1$ & $12$\\\hline
\end{tabular}
\ \ \label{AC1}%
\end{equation}
\textit{or} $q_{3}=3$ and $\xi_{1}+\xi_{2}\in\{0,1,\ldots,9\}$ with $\xi
_{1},\xi_{2}\geq0.$ If in the latter case (\ref{Condition1}) were true, then,
in particular, $q_{1}\mid\widehat{p}_{1}q_{2}+q_{3},$ i.e.,
\begin{align*}
9\xi_{1}+6  &  \mid(6\xi_{1}+5)(9\xi_{2}+6)+3\Longrightarrow3\xi_{1}%
+2\mid(3\xi_{1}+2)(6\xi_{2}+4)+3\xi_{2}+3\\
&  \Longrightarrow3\xi_{1}+2\mid3\xi_{2}+3\Longrightarrow3\xi_{1}+2\mid\xi
_{2}+1,\text{ i.e.,}%
\end{align*}%
\[
(\xi_{1},\xi_{2})\in\left\{
(0,1),(0,3),(0,5),(0,7),(0,9),(1,4),(2,7)\right\}  .
\]
As for everyone of the $7$ possible values of $(\xi_{1},\xi_{2})$ at least one
of the divisibility conditions (\ref{Condition1}) and (\ref{Condition2}) is
violated$,$ there are no further admissible triples of pairs.\smallskip

\noindent{}$\bullet$ Case $\left(  \mathbf{5,3},\mathbf{6}\right)  $: $\xi
_{1}=\xi_{2}=1,$ $q_{3}=9,$ and we obtain just one admissible triple of pairs:%
\begin{equation}%
\begin{tabular}
[c]{|c|c|c|c|c|c|}\hline
$p_{1}$ & $q_{1}$ & $p_{2}$ & $q_{2}$ & $p_{3}$ & $q_{3}$\\\hline\hline
$7$ & $9$ & $4$ & $9$ & $1$ & $9$\\\hline
\end{tabular}
\ \ \ \label{AD1}%
\end{equation}
Working symmetrically with type combinations $(\alpha_{1},\alpha_{2}%
,\alpha_{3}),$ where
\[
\alpha_{1},\alpha_{3}\in\{\mathbf{1},\mathbf{2},\mathbf{3},\mathbf{4}%
,\mathbf{5}\}\ \ \text{and\ \ }\alpha_{2}=\mathbf{6},
\]
we determine the admissible triples of pairs:%
\begin{equation}%
\begin{tabular}
[c]{|c|c|c|c|c|c|}\hline
$p_{1}$ & $q_{1}$ & $p_{2}$ & $q_{2}$ & $p_{3}$ & $q_{3}$\\\hline\hline
$2$ & $3$ & $1$ & $6$ & $2$ & $3$\\\hline
\end{tabular}
\ \ \label{AB2}%
\end{equation}
in the case $(\mathbf{1,6,1}),$
\begin{equation}%
\begin{tabular}
[c]{|c|c|c|c|c|c|}\hline
$p_{1}$ & $q_{1}$ & $p_{2}$ & $q_{2}$ & $p_{3}$ & $q_{3}$\\\hline\hline
$2$ & $3$ & $1$ & $9$ & $5$ & $6$\\\hline
\end{tabular}
\ \ \label{Ex3}%
\end{equation}
in the case $(\mathbf{1,6,4}),$%
\begin{equation}%
\begin{tabular}
[c]{|c|c|c|c|c|c|}\hline
$p_{1}$ & $q_{1}$ & $p_{2}$ & $q_{2}$ & $p_{3}$ & $q_{3}$\\\hline\hline
$4$ & $9$ & $1$ & $9$ & $7$ & $9$\\\hline
\end{tabular}
\ \ \label{AD2}%
\end{equation}
in the case $\left(  \mathbf{3},\mathbf{6,5}\right)  $
\begin{equation}%
\begin{tabular}
[c]{|c|c|c|c|c|c|}\hline
$p_{1}$ & $q_{1}$ & $p_{2}$ & $q_{2}$ & $p_{3}$ & $q_{3}$\\\hline\hline
$5$ & $6$ & $1$ & $9$ & $2$ & $3$\\\hline
\end{tabular}
\ \ \label{Ex4}%
\end{equation}
in the case $(\mathbf{4,6,1}),$ and%
\begin{equation}%
\begin{tabular}
[c]{|c|c|c|c|c|c|}\hline
$p_{1}$ & $q_{1}$ & $p_{2}$ & $q_{2}$ & $p_{3}$ & $q_{3}$\\\hline\hline
$5$ & $6$ & $1$ & $12$ & $5$ & $6$\\\hline
\end{tabular}
\ \ \label{AC2}%
\end{equation}
in the case $(\mathbf{4,6,4}).$\hfill{}$\square$

\section{Proof of Theorem \ref{THM3}: Step 3\label{STEP3}}

\begin{lemma}
\label{LEMMA31}There exist exactly $23$ admissible triples of pairs
\emph{(\ref{TRIPLES}) }among those corresponding to the type combinations
$(\alpha_{1},\alpha_{2},\alpha_{3})$ with $\alpha_{1}\in\{\mathbf{1}%
,\mathbf{2},\mathbf{3},\mathbf{4},\mathbf{5}\}$ and $\alpha_{2},\alpha_{3}%
\in\{\mathbf{6},\mathbf{7}\}.$
\end{lemma}

\begin{proof}
For every $\alpha_{1}\in\{\mathbf{1},\mathbf{2},\mathbf{3},\mathbf{4}%
,\mathbf{5}\}$ we consider the combinations%
\[%
\begin{tabular}
[c]{|c|c|c|c|c|}\hline
\textbf{Case} & $p_{2}$ & $q_{2}$ & $p_{3}=\widehat{p}_{3}$ & $q_{3}%
$\\\hline\hline
$(\alpha_{1}\mathbf{,6,6})$ & $1$ & $\geq2$ & $1$ & $\geq2$\\\hline
$(\alpha_{1}\mathbf{,6,7})$ & $1$ & $\geq2$ & $0$ & $1$\\\hline
$(\alpha_{1}\mathbf{,7,6})$ & $0$ & $1$ & $1$ & $\geq2$\\\hline
$(\alpha_{1}\mathbf{,7,7})$ & $0$ & $1$ & $0$ & $1$\\\hline
\end{tabular}
\ \ \ \
\]
and examine what happens in each of the twenty cases separately.\smallskip
\ \newline$\bullet$ Case $(\mathbf{1,6,6})$: Here, and for the next three
cases, $p_{1}=\widehat{p}_{1}=2,$ $q_{1}=3$ and $s_{1}=1.$ By
(\ref{Condition3}) and (\ref{Condition4}) the pair $\left(  q_{2}%
,q_{3}\right)  $ has to be chosen from the $21$ elements of the set
\[
\left\{  \left.  \left(  q_{2},q_{3}\right)  \in\mathbb{Z}^{2}\right\vert
q_{2}\geq2,\ q_{3}\geq2,\text{ and }q_{2}+q_{3}\leq9\right\}  .
\]
Taking into account the divisibility conditions (\ref{Condition1}),
(\ref{Condition2}), i.e., $3q_{2}\mid2q_{2}+q_{3}+3$ and $3q_{3}\mid
2q_{3}+q_{2}+3,$ we obtain $\left(  q_{2},q_{3}\right)  \in\{(2,5),(5,2)\}.$
Hence, there are two admissible triples of pairs, namely
\begin{equation}%
\begin{tabular}
[c]{|c|c|c|c|c|c|}\hline
$p_{1}$ & $q_{1}$ & $p_{2}$ & $q_{2}$ & $p_{3}$ & $q_{3}$\\\hline\hline
$2$ & $3$ & $1$ & $2$ & $1$ & $5$\\\hline
\end{tabular}
\ \ \ \ \ \ \ \ \label{A1}%
\end{equation}
and
\begin{equation}%
\begin{tabular}
[c]{|c|c|c|c|c|c|}\hline
$p_{1}$ & $q_{1}$ & $p_{2}$ & $q_{2}$ & $p_{3}$ & $q_{3}$\\\hline\hline
$2$ & $3$ & $1$ & $5$ & $1$ & $2$\\\hline
\end{tabular}
\ \ \ \ \ \ \ \ \label{A2}%
\end{equation}
\noindent$\bullet$ Case $(\mathbf{1,6,7})$: By (\ref{Condition3}) (or
(\ref{Condition4})) we have $q_{2}\leq8.$ By (\ref{Condition2}), $3\mid
q_{2}-1,$ i.e., $\ q_{2}\in\{4,7\}.$ The value $q_{2}=7$ does not satisfy
(\ref{Condition1}): $3q_{2}\mid2q_{2}+4.$ Hence, there is only one admissible
triple of pairs:%
\begin{equation}%
\begin{tabular}
[c]{|c|c|c|c|c|c|}\hline
$p_{1}$ & $q_{1}$ & $p_{2}$ & $q_{2}$ & $p_{3}$ & $q_{3}$\\\hline\hline
$2$ & $3$ & $1$ & $4$ & $0$ & $1$\\\hline
\end{tabular}
\ \ \ \ \ \ \ \label{B1}%
\end{equation}
$\bullet$ Case $(\mathbf{1,7,6})$: Analogously, we find just one admissible
triple of pairs:
\begin{equation}%
\begin{tabular}
[c]{|c|c|c|c|c|c|}\hline
$p_{1}$ & $q_{1}$ & $p_{2}$ & $q_{2}$ & $p_{3}$ & $q_{3}$\\\hline\hline
$2$ & $3$ & $0$ & $1$ & $1$ & $4$\\\hline
\end{tabular}
\ \ \ \ \ \ \ \ \ \label{B2}%
\end{equation}
$\bullet$ Case $(\mathbf{1,7,7})$: In this case both divisibility conditions
(\ref{Condition1}) and (\ref{Condition2}) are satisfied automatically and lead
to the admissible triple of pairs:%
\[%
\begin{tabular}
[c]{|c|c|c|c|c|c|}\hline
$p_{1}$ & $q_{1}$ & $p_{2}$ & $q_{2}$ & $p_{3}$ & $q_{3}$\\\hline\hline
$2$ & $3$ & $0$ & $1$ & $0$ & $1$\\\hline
\end{tabular}
\ \ \ \ \ \ \ \ \
\]
$\bullet$ Case $(\mathbf{2,6,6})$: Here, and for the next three cases,
$p_{1}=\widehat{p}_{1}=3\xi_{1}+2,$ $q_{1}=9\xi_{1}+3$ and $s_{1}=\xi_{1}+2$
for an integer $\xi_{1}\geq1.$ By (\ref{Condition3}) we have $\xi_{1}%
+q_{2}+q_{3}\leq9.$ Condition (\ref{Condition1}) reads as%
\[
3(3\xi_{1}+1)q_{2}\mid(3\xi_{1}+2)q_{2}+(9\xi_{1}+3)+q_{3}=3(3\xi_{1}%
+1)q_{2}-6\xi_{1}q_{2}-q_{2}+(9\xi_{1}+3)+q_{3},
\]
i.e.,%
\[
3(3\xi_{1}+1)q_{2}\mid6\xi_{1}q_{2}+q_{2}-(9\xi_{1}+3)-q_{3}\text{ \ with
\ }6\xi_{1}q_{2}+q_{2}-(9\xi_{1}+3)-q_{3}\leq0.
\]
Since $1\leq\xi_{1}\leq5,$
\[
(6\xi_{1}+1)q_{2}\leq9\xi_{1}+3+q_{3}\leq9\xi_{1}+3+\left(  9-\xi_{1}%
-q_{2}\right)  \Longrightarrow(6\xi_{1}+2)q_{2}\leq8\xi_{1}+12\leq52,
\]
implying%
\[
8\leq(3\xi_{1}+1)q_{2}\leq26.
\]
These inequalities are satisfied if and only if%
\[
\left(  \xi_{1},q_{2}\right)  \in\left\{
(1,2),(1,3),(1,4),(1,5),(1,6),(2,2),(2,3),(2,4),(3,2)\right\}  .
\]
Since $2\leq q_{3}\leq9-(\xi_{1}+q_{2}),$ the divisibility condition
(\ref{Condition1}) is true only for $\xi_{1}=1,$ $q_{2}=3$ and $q_{3}=2.$ (For
these values (\ref{Condition2}) is also true.) Hence, the only admissible
triple of pairs is the following:
\[%
\begin{tabular}
[c]{|c|c|c|c|c|c|}\hline
$p_{1}$ & $q_{1}$ & $p_{2}$ & $q_{2}$ & $p_{3}$ & $q_{3}$\\\hline\hline
$5$ & $12$ & $1$ & $2$ & $1$ & $2$\\\hline
\end{tabular}
\ \ \ \ \ \ \ \ \ \ \
\]
\noindent\noindent$\bullet$ Case $(\mathbf{2,6,7})$: By (\ref{Condition3}) we
have $\xi_{1}+q_{2}\leq8.$ Condition (\ref{Condition2}) gives%
\[
3(3\xi_{1}+1)\mid3\xi_{1}+2+q_{2}\Longrightarrow3\mid q_{2}-1\text{ and }%
3\xi_{1}+1\mid q_{2}+1.
\]
But this means that $\left(  \xi_{1},q_{2}\right)  \in\{(1,7),(2,6)\}.$
$(2,6)$ is not permitted because $21\nmid14$ and $(1,7)$ violates
(\ref{Condition1}), so there are no admissible triples of pairs.\smallskip
\newline$\bullet$ Case $(\mathbf{2,7,6})$: As in the case $(\mathbf{2,6,7})$
one shows that there are no admissible triples of pairs.\smallskip
\newline$\bullet$ Case $(\mathbf{2,7,7})$: Conditions (\ref{Condition1}) and
(\ref{Condition2}) give $q_{1}=3(p_{1}-1)\mid p_{1}+1,$ i.e., $p_{1}=2,$ but
in this case $p_{1}\geq5.$ Hence, there are no admissible triples of
pairs.\smallskip\newline$\bullet$ Case $(\mathbf{3,6,6})$: Here, and for the
next three cases, $p_{1}=3\xi_{1}+1,$ $\widehat{p}_{1}=6\xi_{1}+1,$
$q_{1}=9\xi_{1}$ and $s_{1}=\xi_{1}+1$ for an integer $\xi_{1}\geq1.$ By
(\ref{Condition3}) we have $\xi_{1}+q_{2}+q_{3}\leq11.$ Condition
(\ref{Condition1}) reads as%
\[
9\xi_{1}q_{2}\mid(6\xi_{1}+1)q_{2}+9\xi_{1}+q_{3}=9\xi_{1}q_{2}-3\xi_{1}%
q_{2}+q_{2}+9\xi_{1}+q_{3},
\]
i.e.,%
\[
9\xi_{1}q_{2}\mid3\xi_{1}q_{2}-q_{2}-9\xi_{1}-q_{3}\text{ \ with \ }3\xi
_{1}q_{2}-q_{2}-9\xi_{1}-q_{3}\leq0.
\]
Since $1\leq\xi_{1}\leq7,$
\[
3\xi_{1}q_{2}\leq q_{2}+q_{3}+9\xi_{1}\leq11+8\xi_{1}\leq67\Longrightarrow
2\leq\xi_{1}q_{2}\leq22.
\]
Since $2\leq q_{3}\leq11-(\xi_{1}+q_{2}),$ the divisibility conditions
(\ref{Condition1}) and (\ref{Condition2}) are true only for $\xi_{1}=1,$
$q_{2}=6,$ $q_{3}=3,$ or $\xi_{1}=2,$ $q_{2}=4,$ $q_{3}=2,$ leading to two
admissible triple of pairs, namely
\begin{equation}%
\begin{tabular}
[c]{|c|c|c|c|c|c|}\hline
$p_{1}$ & $q_{1}$ & $p_{2}$ & $q_{2}$ & $p_{3}$ & $q_{3}$\\\hline\hline
$4$ & $9$ & $1$ & $6$ & $1$ & $3$\\\hline
\end{tabular}
\ \ \ \ \ \ \ \ \ \label{C1}%
\end{equation}
and%
\begin{equation}%
\begin{tabular}
[c]{|c|c|c|c|c|c|}\hline
$p_{1}$ & $q_{1}$ & $p_{2}$ & $q_{2}$ & $p_{3}$ & $q_{3}$\\\hline\hline
$7$ & $18$ & $1$ & $4$ & $1$ & $2$\\\hline
\end{tabular}
\ \ \ \ \ \ \ \ \ \label{D1}%
\end{equation}
$\bullet$ Case $(\mathbf{3,6,7})$: By (\ref{Condition3}) we have $\xi
_{1}+q_{2}\leq10.$ Condition (\ref{Condition2}) gives%
\[
9\xi_{1}\mid3\xi_{1}+1+q_{2}\Longrightarrow9\xi_{1}\mid6\xi_{1}-q_{2}-1,\text{
with \ }6\xi_{1}-q_{2}-1\leq0.
\]
Thus,
\[6\xi_{1}\leq q_{2}+1\leq11-\xi_{1}\Longrightarrow \xi_{1}\leq
\frac{11}{7}\Longrightarrow \xi_{1}=1.
\]
Since $2\leq q_{2}\leq9,$ condition
(\ref{Condition1}) (i.e., $9q_{2}\mid7q_{2}+10$) implies $q_{2}=5.$ The
corresponding admissible triple of pairs is the following:
\begin{equation}%
\begin{tabular}
[c]{|c|c|c|c|c|c|}\hline
$p_{1}$ & $q_{1}$ & $p_{2}$ & $q_{2}$ & $p_{3}$ & $q_{3}$\\\hline\hline
$4$ & $9$ & $1$ & $5$ & $0$ & $1$\\\hline
\end{tabular}
\ \ \ \ \ \ \ \ \label{F1}%
\end{equation}
$\bullet$ Case $(\mathbf{3,7,6})$: By (\ref{Condition3}) we have $\xi
_{1}+q_{3}\leq10.$ Condition (\ref{Condition1}) gives%
\[
\xi_{1}\mid6\xi_{1}+1+q_{3}\Longrightarrow9\xi_{1}\mid3\xi_{1}-q_{3}-1,\text{
with \ }3\xi_{1}-q_{3}-1\leq0.
\]
Thus,
\[3\xi_{1}\leq q_{3}+1\leq11-\xi_{1}\Longrightarrow \xi_{1}\leq
\frac{11}{4}\Longrightarrow \xi_{1}\in\{1,2\}.
\]
Since $2\leq q_{3}\leq9,$
condition (\ref{Condition1}) implies $\left(  \xi_{1},q_{3}\right)
\in\left\{  \left(  1,2\right)  ,(2,5)\right\}  .$ $(2,5)$ is not permitted
because it violates (\ref{Condition2}). For this reason, the only admissible
triple of pairs is the following:
\begin{equation}%
\begin{tabular}
[c]{|c|c|c|c|c|c|}\hline
$p_{1}$ & $q_{1}$ & $p_{2}$ & $q_{2}$ & $p_{3}$ & $q_{3}$\\\hline\hline
$4$ & $9$ & $0$ & $1$ & $1$ & $2$\\\hline
\end{tabular}
\ \ \ \ \ \ \ \ \ \ \label{E1}%
\end{equation}
$\bullet$ Case $(\mathbf{3,7,7})$: Condition (\ref{Condition2}) gives
$q_{1}=3(p_{1}-1)\mid p_{1}+1,$ i.e., $p_{1}=2,$ but in this case $p_{1}%
\geq4.$ Hence, there are no admissible triples of pairs.\smallskip
\newline$\bullet$ Case $(\mathbf{4,6,6})$: Here, and for the next three cases,
$p_{1}=\widehat{p}_{1}=6\xi_{1}+5,$ $q_{1}=9\xi_{1}+6$ and $s_{1}=\xi_{1}+1$
for an integer $\xi_{1}\geq0.$ By (\ref{Condition3}) we have $\xi_{1}%
+q_{2}+q_{3}\leq11.$ Condition (\ref{Condition1}) reads as%
\[
(9\xi_{1}+6)q_{2}\mid(6\xi_{1}+5)q_{2}+9\xi_{1}+6+q_{3}=(9\xi_{1}%
+6)q_{2}-(3\xi_{1}+1)q_{2}+9\xi_{1}+6+q_{3},
\]
i.e.,%
\[
(9\xi_{1}+6)q_{2}\mid(3\xi_{1}+1)q_{2}-9\xi_{1}-6-q_{3}\text{ \ with \ }%
(3\xi_{1}+1)q_{2}-9\xi_{1}-6-q_{3}\leq0.
\]
Since $1\leq\xi_{1}\leq7,$ we obtain
\[
(3\xi_{1}+1)q_{2}\leq9\xi_{1}+6+(11-q_{2}-\xi_{1})\Longrightarrow(3\xi
_{1}+2)q_{2}\leq17+8\xi_{1}\leq73,
\]
i.e., $4\leq(3\xi_{1}+2)q_{2}\leq73.$ Since $2\leq q_{3}\leq11-(\xi_{1}%
+q_{2}),$ the divisibility conditions (\ref{Condition1}) and (\ref{Condition2}%
) are satisfied only for $\xi_{1}\in\{0,1,2\}.$ In particular$,$ for $\xi
_{1}=0$ we obtain $\left(  q_{2},q_{3}\right)  \in\left\{  \left(  8,2\right)
,(2,8)\right\}  $ and the admissible triples of pairs
\begin{equation}%
\begin{tabular}
[c]{|c|c|c|c|c|c|}\hline
$p_{1}$ & $q_{1}$ & $p_{2}$ & $q_{2}$ & $p_{3}$ & $q_{3}$\\\hline\hline
$5$ & $6$ & $1$ & $8$ & $1$ & $2$\\\hline
\end{tabular}
\ \ \ \ \ \ \ \ \ \label{G1}%
\end{equation}
and%
\begin{equation}%
\begin{tabular}
[c]{|c|c|c|c|c|c|}\hline
$p_{1}$ & $q_{1}$ & $p_{2}$ & $q_{2}$ & $p_{3}$ & $q_{3}$\\\hline\hline
$5$ & $6$ & $1$ & $2$ & $1$ & $8$\\\hline
\end{tabular}
\ \ \ \ \ \ \ \ \ \label{G2}%
\end{equation}
For $\xi_{1}=1$ we have necessarily $q_{2}=q_{3}=5$ and the admissible triple
of pairs:
\[%
\begin{tabular}
[c]{|c|c|c|c|c|c|}\hline
$p_{1}$ & $q_{1}$ & $p_{2}$ & $q_{2}$ & $p_{3}$ & $q_{3}$\\\hline\hline
$11$ & $15$ & $1$ & $5$ & $1$ & $5$\\\hline
\end{tabular}
\ \ \ \ \ \ \ \ \
\]
Finally, for $\xi_{1}=2$ we have necessarily $q_{2}=q_{3}=4$ and the
admissible triple of pairs:
\[%
\begin{tabular}
[c]{|c|c|c|c|c|c|}\hline
$p_{1}$ & $q_{1}$ & $p_{2}$ & $q_{2}$ & $p_{3}$ & $q_{3}$\\\hline\hline
$17$ & $24$ & $1$ & $4$ & $1$ & $4$\\\hline
\end{tabular}
\ \ \ \ \ \ \ \ \
\]
$\bullet$ Case $(\mathbf{4,6,7})$: By (\ref{Condition3}) we have $\xi
_{1}+q_{2}\leq10.$ Condition (\ref{Condition2}) gives%
\[
9\xi_{1}+6\mid6\xi_{1}+5+q_{2}\Longrightarrow9\xi_{1}+6\mid3\xi_{1}%
+1-q_{2},\text{ with \ }3\xi_{1}+1-q_{2}\leq0.
\]
Thus, $3\xi_{1}\leq q_{2}-1\leq9-\xi_{1}\Longrightarrow$ $\xi_{1}\leq\frac
{9}{4}\Longrightarrow$ $\xi_{1}\in\{0,1,2\}.$ Since $2\leq q_{2}\leq10,$
condition (\ref{Condition2}) implies
\[
(\xi_{1},q_{2})\in\left\{  \left(  0,7\right)  ,\left(  1,4\right)  ,\left(
2,7\right)  \right\}  .
\]
$\left(  2,7\right)  $ is not permitted because it violates (\ref{Condition1}%
); therefore, the admissible triples of pairs are%
\begin{equation}%
\begin{tabular}
[c]{|c|c|c|c|c|c|}\hline
$p_{1}$ & $q_{1}$ & $p_{2}$ & $q_{2}$ & $p_{3}$ & $q_{3}$\\\hline\hline
$5$ & $6$ & $1$ & $7$ & $0$ & $1$\\\hline
\end{tabular}
\ \ \ \ \ \ \ \ \ \label{H1}%
\end{equation}
and%
\begin{equation}%
\begin{tabular}
[c]{|c|c|c|c|c|c|}\hline
$p_{1}$ & $q_{1}$ & $p_{2}$ & $q_{2}$ & $p_{3}$ & $q_{3}$\\\hline\hline
$11$ & $15$ & $1$ & $4$ & $0$ & $1$\\\hline
\end{tabular}
\ \ \ \ \ \ \ \ \label{I1}%
\end{equation}
$\bullet$ Case $(\mathbf{4,7,6})$: As in the case $(\mathbf{4,6,7})$ one
proves that there are two admissible triples of pairs, namely%
\begin{equation}%
\begin{tabular}
[c]{|c|c|c|c|c|c|}\hline
$p_{1}$ & $q_{1}$ & $p_{2}$ & $q_{2}$ & $p_{3}$ & $q_{3}$\\\hline\hline
$5$ & $6$ & $0$ & $1$ & $1$ & $7$\\\hline
\end{tabular}
\ \ \ \ \ \ \ \ \ \label{H2}%
\end{equation}
and%
\begin{equation}%
\begin{tabular}
[c]{|c|c|c|c|c|c|}\hline
$p_{1}$ & $q_{1}$ & $p_{2}$ & $q_{2}$ & $p_{3}$ & $q_{3}$\\\hline\hline
$11$ & $15$ & $0$ & $1$ & $1$ & $4$\\\hline
\end{tabular}
\ \ \ \ \ \ \ \ \label{I2}%
\end{equation}
$\bullet$ Case $(\mathbf{4,7,7})$: Conditions (\ref{Condition1}) and
(\ref{Condition2}) give $q_{1}=\frac{3}{2}(p_{1}-1)\mid p_{1}+1,$ i.e.,
$p_{1}=5.$ Thus, we find just one admissible triple of pairs:
\[%
\begin{tabular}
[c]{|c|c|c|c|c|c|}\hline
$p_{1}$ & $q_{1}$ & $p_{2}$ & $q_{2}$ & $p_{3}$ & $q_{3}$\\\hline\hline
$5$ & $6$ & $0$ & $1$ & $0$ & $1$\\\hline
\end{tabular}
\ \ \ \ \ \ \ \ \
\]
$\bullet$ Case $(\mathbf{5,6,6})$: Here, and for the next three cases,
$p_{1}=6\xi_{1}+1,$ $\widehat{p}_{1}=3\xi_{1}+1,$ $q_{1}=9\xi_{1}$ and
$s_{1}=\xi_{1}+1$ for an integer $\xi_{1}\geq1.$ By (\ref{Condition3}) we have
$\xi_{1}+q_{2}+q_{3}\leq11.$ Condition (\ref{Condition1}) reads as%
\[
9\xi_{1}q_{2}\mid(3\xi_{1}+1)q_{2}+9\xi_{1}+q_{3}=9\xi_{1}q_{2}-6\xi_{1}%
q_{2}+q_{2}+9\xi_{1}+q_{3},
\]
i.e.,%
\[
9\xi_{1}q_{2}\mid6\xi_{1}q_{2}-q_{2}-9\xi_{1}-q_{3}\text{ \ with \ }6\xi
_{1}q_{2}-q_{2}-9\xi_{1}-q_{3}\leq0.
\]
Since $1\leq\xi_{1}\leq7,$
\[
6\xi_{1}q_{2}\leq q_{2}+q_{3}+9\xi_{1}\leq11+8\xi_{1}\leq67\Longrightarrow
2\leq\xi_{1}q_{2}\leq11.
\]
Since $2\leq q_{3}\leq11-(\xi_{1}+q_{2}),$ the divisibility conditions
(\ref{Condition1}) and (\ref{Condition2}) are satisfied only for $\xi_{1}=1,$
$q_{2}=3,$ $q_{3}=6,$ or $\xi_{1}=2,$ $q_{2}=2,$ $q_{3}=4,$ leading to two
admissible triple of pairs, namely
\begin{equation}%
\begin{tabular}
[c]{|c|c|c|c|c|c|}\hline
$p_{1}$ & $q_{1}$ & $p_{2}$ & $q_{2}$ & $p_{3}$ & $q_{3}$\\\hline\hline
$7$ & $9$ & $1$ & $3$ & $1$ & $6$\\\hline
\end{tabular}
\ \ \ \ \ \ \ \ \ \ \ \label{C2}%
\end{equation}
and
\begin{equation}%
\begin{tabular}
[c]{|c|c|c|c|c|c|}\hline
$p_{1}$ & $q_{1}$ & $p_{2}$ & $q_{2}$ & $p_{3}$ & $q_{3}$\\\hline\hline
$13$ & $18$ & $1$ & $2$ & $1$ & $4$\\\hline
\end{tabular}
\ \ \ \ \ \ \ \ \ \ \label{D2}%
\end{equation}
$\bullet$ Case $(\mathbf{5,6,7})$: By (\ref{Condition3}) we have $\xi
_{1}+q_{2}\leq10.$ Condition (\ref{Condition2}) gives%
\[
9\xi_{1}\mid6\xi_{1}+1+q_{2}\Longrightarrow9\xi_{1}\mid3\xi_{1}-q_{2}-1,\text{
with \ }3\xi_{1}-q_{2}-1\leq0.
\]
Thus,
\[
3\xi_{1}\leq q_{2}+1\leq11-\xi_{1}\Longrightarrow \xi_{1}\leq
\frac{11}{4}\Longrightarrow \xi_{1}\in\{1,2\}.
\]
Since $2\leq q_{2}\leq9,$
condition (\ref{Condition2}) implies $\xi_{1}=1$ and $q_{2}=2.$ The result is
the following admissible triple of pairs:%
\begin{equation}%
\begin{tabular}
[c]{|c|c|c|c|c|c|}\hline
$p_{1}$ & $q_{1}$ & $p_{2}$ & $q_{2}$ & $p_{3}$ & $q_{3}$\\\hline\hline
$7$ & $9$ & $1$ & $2$ & $0$ & $1$\\\hline
\end{tabular}
\ \ \ \ \ \ \ \ \label{E2}%
\end{equation}
$\bullet$ Case $(\mathbf{5,7,6})$: By (\ref{Condition3}), $\xi_{1}+q_{3}%
\leq10.$ Now (\ref{Condition1}) reads as%
\[
9\xi_{1}\mid3\xi_{1}+1+q_{3}\Longrightarrow9\xi_{1}\mid6\xi_{1}-q_{3}-1,\text{
with \ }6\xi_{1}-q_{3}-1\leq0.
\]
Thus,
\[6\xi_{1}\leq q_{3}+1\leq11-\xi_{1}\Longrightarrow \xi_{1}\leq
\frac{11}{7}\Longrightarrow \xi_{1}=1.
\]
Since $2\leq q_{3}\leq9,$ condition
(\ref{Condition1}) implies $q_{3}=5.$ The corresponding admissible triple of
pairs is the following:
\begin{equation}%
\begin{tabular}
[c]{|c|c|c|c|c|c|}\hline
$p_{1}$ & $q_{1}$ & $p_{2}$ & $q_{2}$ & $p_{3}$ & $q_{3}$\\\hline\hline
$7$ & $9$ & $0$ & $1$ & $1$ & $5$\\\hline
\end{tabular}
\ \ \ \ \ \ \ \ \label{F2}%
\end{equation}
$\bullet$ Case $(\mathbf{5,7,7})$: Condition (\ref{Condition2}) gives
$q_{1}=\frac{3}{2}(p_{1}-1)\mid p_{1}+1,$ i.e., $p_{1}\leq5,$ but in this case
$p_{1}\geq7.$ Hence, there are no admissible triples of pairs.
\end{proof}

\begin{remark}
The majority of the admissible triples of pairs induce toric log Del Pezzo
surfaces admitting at least one Gorenstein singularity. This is due to the
fact that the $q_{i}$'s corresponding to Gorenstein singularities can be
viewed as parameters moving freely between $2$ and an upper bound dictated by
conditions (\ref{Condition3}) and (\ref{Condition4}), without any further restrictions.
\end{remark}

\section{Proof of Theorem \ref{THM3}: Step4\label{STEP4}}

\begin{lemma}
\label{ISOLEMMA}The toric log Del Pezzo surfaces induced by the following
admissible triples of pairs $\mathbf{(a)}$ and $\mathbf{(b)}$\emph{:}
\[%
\begin{tabular}
[c]{|c|c|c|c|c|c|c|}\hline
$\mathbf{(a)}$ & \emph{(\ref{AB1})} & \emph{(\ref{AC1})} & \emph{(\ref{AD1})}
& \emph{(\ref{A1})} & \emph{(\ref{B1})} & \emph{(\ref{C1})}\\\hline
$\mathbf{(b)}$ & \emph{(\ref{AB2})} & \emph{(\ref{AC2})} & \emph{(\ref{AD2})}
& \emph{(\ref{A2})} & \emph{(\ref{B2})} & \emph{(\ref{C2})}\\\hline\hline
$\mathbf{(a)}$ & \emph{(\ref{D1})} & \emph{(\ref{E1})} & \emph{(\ref{F1})} &
\emph{(\ref{G1})} & \emph{(\ref{H1})} & \emph{(\ref{I1})}\\\hline
$\mathbf{(b)}$ & \emph{(\ref{D2})} & \emph{(\ref{E2})} & \emph{(\ref{F2})} &
\emph{(\ref{G2})} & \emph{(\ref{H2})} & \emph{(\ref{I2})}\\\hline
\end{tabular}
\
\]
are isomorphic to each other. The same is true for the four surfaces induced
by the following admissible triples of pairs\emph{:}
\[%
\begin{tabular}
[c]{|c|c|c|c|}\hline
$\mathbf{(a)}$ & $\mathbf{(b)}$ & $\mathbf{(c)}$ & $\mathbf{(d)}%
$\\\hline\hline
\emph{(\ref{Ex1})} & \emph{(\ref{Ex2})} & \emph{(\ref{Ex3})} &
\emph{(\ref{Ex4})}\\\hline
\end{tabular}
\]
\emph{(The admissible triples of pairs are given by their reference numbers.)}
\end{lemma}

\begin{proof}
If $X_{\Delta_{\mathbf{(a)}}}$ (resp., $X_{\Delta_{\mathbf{(b)}}}$) is the
toric Del Pezzo surface induced by the admissible triple of pairs
$\mathbf{(a)}$ (resp., $\mathbf{(b)}$) in the first list, then $\mathfrak{G}%
_{\Delta_{\mathbf{(a)}}}\overset{\text{gr.}}{\cong}\mathfrak{G}_{\Delta
_{\mathbf{(b)}}}^{\text{rev}}.$ Correspondingly, if $X_{\Delta_{\mathbf{(a)}}%
},$ $X_{\Delta_{\mathbf{(b)}}},$ $X_{\Delta_{\mathbf{(c)}}},$ $X_{\Delta
_{\mathbf{(d)}}}$ are the four surfaces induced by the admissible triples of
pairs in the second list, then we obtain
\[
\mathfrak{G}_{\Delta_{\mathbf{(a)}}}\overset{\text{gr.}}{\cong}\mathfrak{G}%
_{\Delta_{\mathbf{(b)}}}^{\text{rev}}\overset{\text{gr.}}{\cong}%
\mathfrak{G}_{\Delta_{\mathbf{(c)}}}^{\text{rev}}\overset{\text{gr.}}{\cong%
}\mathfrak{G}_{\Delta_{\mathbf{(d)}}}.
\]
It is therefore enough to apply Theorem \ref{CLASSIFTHM}.
\end{proof}

\begin{note}
By Lemmas \ref{Lemma12}, \ref{Lemma21}, \ref{Lemma22} and \ref{LEMMA31} we
proved that among all possible triples of pairs there exist exactly $33$ which are
admissible. Lemma \ref{ISOLEMMA} informs us that, in fact, for the
classification of toric Del Pezzo surfaces $X_{\Delta}$ having Picard number
$\rho\left(  X_{\Delta}\right)  =1$ and index $\ell=3$ \textit{up to
isomorphism}, we need only $18$ out of them. (The $X_{\Delta}$'s induced by
such a choice of $18$ admissible triples of pairs are \textit{obviously}
pairwise non-isomorphic.)
\end{note}

\noindent{}\textit{End of the proof of Theorem }\ref{THM3}: We consider $18$
representatives of admissible triple of pairs inducing pairwise non-isomorphic
toric Del Pezzo surfaces $X_{\Delta}$ with $\rho\left(  X_{\Delta}\right)  =1$
and index $\ell=3,$ and we enumerate them, e.g., as in the Table \ref{Table6}.
The coordinates of the third minimal generator $\mathbf{n}_{3}$ is computed by
(\ref{N3}). The integers $r_{i}=-\,\overline{C}_{i}^{2},i\in\{1,2,3\},$ are
computed directly via (\ref{CONDri}).

\begin{figure}[h!]
\epsfig{file=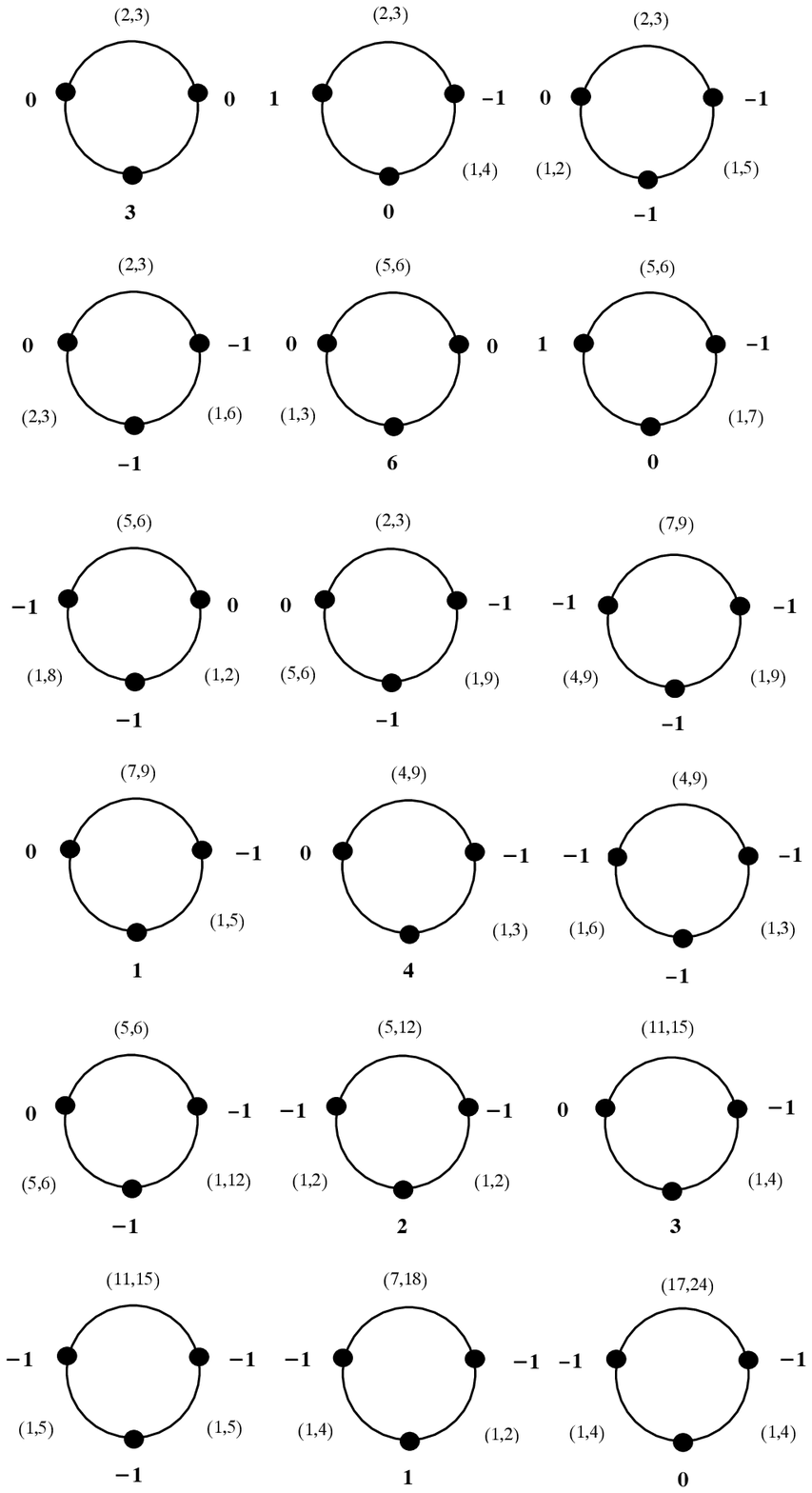, height=19cm, width=13cm}
\caption{}\label{Fig.3}
\end{figure}

\begin{table}[h]
{\small
\[%
\begin{tabular}
[c]{|c|c|c|c|c|c|c|c|c|c|c|c|}\hline
\textbf{No.} & \textbf{Case} & $p_{1}$ & $q_{1}$ & $p_{2}$ & $q_{2}$ & $p_{3}$
& $q_{3}$ & $\mathbf{n}_{3}$ & $r_{1}$ & $r_{2}$ & $r_{3}$\\\hline\hline
(i) & $(\mathbf{1,7,7})$ & $2$ & $3$ & $0$ & $1$ & $0$ & $1$ & $(-1,-1)$ & $0$
& $0$ & $-3$\\\hline
(ii) & $(\mathbf{1,7,6})$ & $2$ & $3$ & $0$ & $1$ & $1$ & $4$ & $(-3,-4)$ &
$1$ & $-1$ & $0$\\\hline
(iii) & $(\mathbf{1,6,6})$ & $2$ & $3$ & $1$ & $2$ & $1$ & $5$ & $(-4,-5)$ &
$1$ & $0$ & $1$\\\hline
(iv) & $(\mathbf{1,1,6})$ & $2$ & $3$ & $2$ & $3$ & $1$ & $6$ & $(-5,-6)$ &
$1$ & $0$ & $1$\\\hline
(v) & $(\mathbf{4,7,7})$ & $5$ & $6$ & $0$ & $1$ & $0$ & $1$ & $(-1,-1)$ & $0$
& $0$ & $-6$\\\hline
(vi) & $(\mathbf{4,7,6})$ & $5$ & $6$ & $0$ & $1$ & $1$ & $7$ & $(-6,-7)$ &
$1$ & $-1$ & $0$\\\hline
(vii) & $(\mathbf{4,6,6})$ & $5$ & $6$ & $1$ & $8$ & $1$ & $2$ & $(-3,-2)$ &
$0$ & $1$ & $1$\\\hline
(viii) & $(\mathbf{1,4,6})$ & $2$ & $3$ & $5$ & $6$ & $1$ & $9$ & $(-8,-9)$ &
$1$ & $0$ & $1$\\\hline
(ix) & $(\mathbf{5,3,6})$ & $7$ & $9$ & $4$ & $9$ & $1$ & $9$ & $(-8,-9)$ &
$1$ & $1$ & $1$\\\hline
(x) & $(\mathbf{5,7,6})$ & $7$ & $9$ & $0$ & $1$ & $1$ & $5$ & $(-4,-5)$ & $1$
& $0$ & $-1$\\\hline
(xi) & $(\mathbf{3,7,6})$ & $4$ & $9$ & $0$ & $1$ & $1$ & $2$ & $(-1,-2)$ &
$1$ & $0$ & $-4$\\\hline
(xii) & $(\mathbf{3,6,6})$ & $4$ & $9$ & $1$ & $6$ & $1$ & $3$ & $(-2,-3)$ &
$1$ & $1$ & $1$\\\hline
(xiii) & $(\mathbf{4,4,6})$ & $5$ & $6$ & $5$ & $6$ & $1$ & $12$ & $(-11,-12)$
& $1$ & $0$ & $1$\\\hline
(xiv) & $(\mathbf{2,6,6})$ & $5$ & $12$ & $1$ & $2$ & $1$ & $2$ & $(-1,-2)$ &
$1$ & $1$ & $-2$\\\hline
(xv) & $(\mathbf{4,7,6})$ & $11$ & $15$ & $0$ & $1$ & $1$ & $4$ & $(-3,-4)$ &
$1$ & $0$ & $-3$\\\hline
(xvi) & $(\mathbf{4,6,6})$ & $11$ & $15$ & $1$ & $5$ & $1$ & $5$ & $(-4,-5)$ &
$1$ & $1$ & $1$\\\hline
(xvii) & $(\mathbf{3,6,6})$ & $7$ & $18$ & $1$ & $4$ & $1$ & $2$ & $(-1,-2)$ &
$1$ & $1$ & $-1$\\\hline
(xviii) & $(\mathbf{4,6,6})$ & $17$ & $24$ & $1$ & $4$ & $1$ & $4$ & $(-3,-4)$
& $1$ & $1$ & $0$\\\hline
\end{tabular}
\ \ \ \ \ \
\]
}\smallskip\caption{{}}%
\label{Table6}%
\end{table}\noindent The \textsc{wve}$^{2}$\textsc{c}-graphs $\mathfrak{G}%
_{\Delta}$ (associated to the $18$ $\Delta$'s) are depicted in Figure
\ref{Fig.3} in this order. (The reference to the double weight $(0,1)$ at an
edge of $\mathfrak{G}_{\Delta}$ is always omitted.) Finally, we may identify
the corresponding $X_{\Delta}$'s with weighted projective planes or quotients
thereof by a finite abelian group $H_{\Delta}$ via Lemma \ref{WEIGHTEDLEMMA}.
(In the statement of the Theorem we have w.l.o.g. rearranged the weights in
ascending order. Computing the Smith normal form, $H_{\Delta}$ turns out to be cyclic for the
surfaces (ix) and (xviii)). \hfill{}$\square$

\newpage

\newpage

\begin{thebibliography}{99}                                                                                               %


\bibitem {A-Br}\textsc{Alexeev V. \& Brion M.}: \textit{Boundedness of
spherical Fano varieties}, Proc. of the \textquotedblleft Fano
Conference\textquotedblright, Turin, (2004), pp. 69-80.

\bibitem {Borisov}\textsc{Borisov A.}: \textit{Boundedness of Fano threefolds
with log-terminal singularities of given index}, J. Math. Sci. Univ. Tokyo
\textbf{8} (2001), 329-342.

\bibitem {Bo-Bo}\textsc{Borisov A. \& Borisov L.}: \textit{Singular toric Fano
varieties}, Izvestija Acad. Sci. USSR Sb. Math. \textbf{75} (1993), 227-283.

\bibitem {Conrads}\textsc{Conrads H.}: \textit{Weighted projective spaces and
reflexive simplices}, Manuscripta Math. \textbf{107} (2002), 215-227.

\bibitem {Dais}\textsc{Dais D.I.}: \textit{Geometric combinatorics in the
study of compact toric surfaces. }In \textquotedblleft Algebraic and Geometric
Combinatorics\textquotedblright\ (edited by C. Athanasiadis et. al.),
Contemporary Mathematics, Vol. \textbf{423}, American Mathematical Society,
2007, pp. 71-123.

\bibitem {DN}\textsc{Dais D.I. \& Nill B.}: \textit{A boundedness result for
toric log Del Pezzo surfaces,} archiv:math.AG/ 0707.4567, preprint, 2007.

\bibitem {Ewald}\textsc{Ewald G.}: \textit{Combinatorial Convexity and
Algebraic Geometry}, Graduate Texts in Mathematics, Vol. \textbf{168},
Springer-Verlag, 1996.

\bibitem {Fulton}\textsc{Fulton W.}: \textit{Introduction to Toric Varieties},
Annals of Mathematics Studies, Vol. \textbf{131}, Princeton University Press, 1993.

\bibitem {Hensley}\textsc{Hensley D.}: \textit{Lattice vertex polytopes with
interior lattice points}, Pacific Jour. of Math. \textbf{105} (1983), 183-191.

\bibitem {Hirzebruch1}\textsc{Hirzebruch F.}: \textit{\"{U}ber eine Klasse von
einfach-zusammenh\"{a}ngenden komplexen Mannigfaltigkeiten}, Math. Ann.
\textbf{124}, (1951), 1-22. [See also: \textquotedblleft Gesammelte
Abhandlungen\textquotedblright, Band \textbf{I}, Springer-Verlag, 1987, pp. 1-11.]

\bibitem {Hirzebruch2}\bysame, \textit{\"{U}ber vierdimensionale Riemannsche
Fl\"{a}chen mehrdeutiger analytischer Funktionen von zwei komplexen
Ver\"{a}nderlichen}, Math. Ann. \textbf{126}, (1953), 1-22. [See also:
\textquotedblleft Gesammelte Abhandlungen\textquotedblright, Band \textbf{I},
Springer-Verlag, 1987, pp. 11-32.]

\bibitem {Lag-Zieg}\textsc{Lagarias J. \& Ziegler G.M.}: \textit{Bounds for
lattice polytopes containing a fixed number of interior points in a
sublattice}, Canadian J. Math. \textbf{43} (1991), 1022-1035.

\bibitem {Oda}\textsc{Oda T.}:\textit{\ \ Convex Bodies and Algebraic
Geometry. An Introduction to the Theory of Toric Varieties}. Erg. der Math.
und ihrer Grenzgebiete, 3 Folge, Bd. \textbf{15}, Springer-Verlag, 1988.

\bibitem {Sakai1}\textsc{Sakai F.}: \textit{Anticanonical models of rational
surfaces}, Math. Ann. \textbf{269} (1984), 389-410.

\bibitem {Sakai2}\bysame, \textit{Weil divisors on normal surfaces}, Duke
Math. Jour. \textbf{51} (1984), 877-887.

\bibitem {Sakai3}\bysame, \textit{The structure of normal surfaces}, Duke
Math. Jour. \textbf{52} (1985), 627-648.

\bibitem {Scott}\textsc{Scott P.R.}: \textit{On convex lattice polygons},
Bull. Austral. Math. Soc. \textbf{15} (1976), 395-399.
\end{thebibliography}
\end{document}